\documentclass[11pt]{amsart}
\begin{document}
\theoremstyle{plain}
\newtheorem{thm}{Theorem}[section]
\newtheorem*{thm*}{Theorem}
\newtheorem{prop}[thm]{Proposition}
\newtheorem*{prop*}{Proposition}
\newtheorem{lemma}[thm]{Lemma}
\newtheorem{cor}[thm]{Corollary}
\newtheorem*{conj*}{Conjecture}
\newtheorem*{cor*}{Corollary}
\newtheorem{defn}[thm]{Definition}
\theoremstyle{definition}
\newtheorem*{defn*}{Definition}
\newtheorem{rems}[thm]{Remarks}
\newtheorem*{rems*}{Remarks}
\newtheorem*{proof*}{Proof}
\newtheorem*{not*}{Notation}
\newcommand{\npartial}{\slash\!\!\!\partial}
\newcommand{\Heis}{\operatorname{Heis}}
\newcommand{\Solv}{\operatorname{Solv}}
\newcommand{\Spin}{\operatorname{Spin}}
\newcommand{\SO}{\operatorname{SO}}
\newcommand{\ind}{\operatorname{ind}}
\newcommand{\Index}{\operatorname{index}}
\newcommand{\ch}{\operatorname{ch}}
\newcommand{\rank}{\operatorname{rank}}

\newcommand{\abs}[1]{\lvert#1\rvert}

\newcommand{\blankbox}[2]{%
  \parbox{\columnwidth}{\centering
    \setlength{\fboxsep}{0pt}%
    \fbox{\raisebox{0pt}[#2]{\hspace{#1}}}%
  }%
}
\newcommand{\supp}[1]{\operatorname{#1}}
\newcommand{\norm}[1]{\parallel\, #1\, \parallel}
\newcommand{\ip}[2]{\langle #1,#2\rangle}
\setlength{\parskip}{.3cm}
\newcommand{\nc}{\newcommand}
\nc{\nt}{\newtheorem}
\nc{\gf}[2]{\genfrac{}{}{0pt}{}{#1}{#2}}
\nc{\mb}[1]{{\mbox{$ #1 $}}}
\nc{\real}{{\mathbb R}}
\nc{\comp}{{\mathbb C}}
\nc{\ints}{{\mathbb Z}}
\nc{\Ltoo}{\mb{L^2({\mathbf H})}}
\nc{\rtoo}{\mb{{\mathbf R}^2}}
\nc{\slr}{{\mathbf {SL}}(2,\real)}
\nc{\slz}{{\mathbf {SL}}(2,\ints)}
\nc{\su}{{\mathbf {SU}}(1,1)}
\nc{\so}{{\mathbf {SO}}}
\nc{\hyp}{{\mathbb H}}
\nc{\disc}{{\mathbf D}}
\nc{\torus}{{\mathbb T}}
\newcommand{\tk}{\widetilde{K}}
\newcommand{\boe}{{\bf e}}\newcommand{\bt}{{\bf t}}
\newcommand{\vth}{\vartheta}
\newcommand{\CGh}{\widetilde{\CG}}
\newcommand{\db}{\overline{\partial}}
\newcommand{\tE}{\widetilde{E}}
\newcommand{\tr}{\mbox{tr}}
\newcommand{\ta}{\widetilde{\alpha}}
\newcommand{\tb}{\widetilde{\beta}}
\newcommand{\txi}{\widetilde{\xi}}
\newcommand{\hV}{\hat{V}}
\newcommand{\IC}{\mathbf{C}}
\newcommand{\IZ}{\mathbf{Z}}
\newcommand{\IP}{\mathbf{P}}
\newcommand{\IN}{\mathbf{N}}
\newcommand{\IR}{\mathbf{R}}
\newcommand{\IH}{\mathbf{H}}
\newcommand{\IG}{\mathbf{G}}
\newcommand{\CC}{{\mathcal C}}
\newcommand{\CD}{{\mathcal D}}
\newcommand{\CS}{{\mathcal S}}
\newcommand{\CG}{{\mathcal G}}
\newcommand{\CL}{{\mathcal L}}
\newcommand{\CO}{{\mathcal O}}
\nc{\ca}{{\mathcal A}}
\nc{\cag}{{{\mathcal A}^\Gamma}}
\nc{\cg}{{\mathcal G}}
\nc{\chh}{{\mathcal H}}
\nc{\ck}{{\mathcal B}}
\nc{\cl}{{\mathcal L}}
\nc{\cm}{{\mathcal M}}
\nc{\cn}{{\mathcal N}}
\nc{\cs}{{\mathcal S}}
\nc{\cz}{{\mathcal Z}}
\nc{\sind}{\sigma{\rm -ind}}


\title{}

\titlepage

\begin{center}{\bf SPECTRAL FLOW AND DIXMIER TRACES}
\\
\vspace{.25 in}
{\bf by}\\
\vspace{.25 in}
{\bf Alan Carey}$^{1,3}$\\Department of Pure Mathematics\\University of Adelaide\\
Adelaide, S.A. 5005\\AUSTRALIA\\
\vspace{.25 in}
{\bf John Phillips}$^2$\\Department of Mathematics and Statistics\\
University of Victoria\\Victoria, B.C. V8W 3P4\\CANADA\\
\vspace{.25 in}
{\bf Fyodor Sukochev}$^1$\\School of Informatics and Engineering\\
Flinders University\\Bedford Park 5042\\AUSTRALIA\\
\vspace{0.5 in}
AMS classification nos.
 Primary:
19K56, 46L80: secondary: 58B30, 46L87.

\vspace{0.25 in}
Supported by grants from ARC$^1$ (Australia), the New Zealand
Mathematics Research Institute$^1$ and
NSERC$^2$ (Canada). Part of this research was completed for the
Clay Mathematics Institute$^3$

\end{center}
\begin{abstract}
We obtain general theorems which enable the calculation of the
Dixmier trace in terms of the asymptotics of the zeta function
and of the heat operator in a general semi-finite von Neumann
algebra. Our results have several applications.  We deduce a formula for
the Chern character of
an odd ${\mathcal L}^{(1,\infty)}$-summable Breuer-Fredholm module
in terms of a Hochschild 1-cycle.
We explain how to derive
a Wodzicki residue for pseudo-differential operators along the
orbits of an ergodic $\IR^n$ action on a compact space $X$. Finally
we give a short proof an index theorem of Lesch for generalised
Toeplitz operators.
\end{abstract}
\maketitle
\newpage

\section{\bf Introduction}

There is a generalisation of the usual setting of noncommutative
geometry where one replaces spectral triples
by Breuer-Fredholm modules. In this situation one is given a Hilbert space
$\mathcal H$, a $C^*$-algebra $\mathcal A$ represented
in a semifinite von Neumann algebra $\mathcal N$
which acts on $\mathcal H$ and a self-adjoint unbounded
operator $D_0$ affiliated to  $\mathcal N$ and such
that the commutator $[a,D_0]$ is bounded for a dense
set of $a\in \mathcal A$ \cite{CPS}. This situation arises
for example in the twisted $L^2$-index theorem of Gromov \cite{Gr}.
There are also other interesting
invariants of operators affiliated to $\mathcal N$
such as $L^2$ spectral flow studied in \cite{CP1}, \cite{CP2}.
We became interested in the Dixmier trace
and its relation to the zeta function partly as
a result of the local index formula of
Connes and Moscovici \cite{CM}. In \cite{CM}
a formula for spectral flow in an ${\mathcal L}^{(p,\infty)}$-summable
Fredholm module (the notation for these symmetric ideals is
explained below) is given.
It is natural to try to relate this formula and those
for spectral flow in \cite{CP1}, \cite{CP2}.

In the course of this investigation we became aware of
the subtleties in the zeta function approach to the Dixmier trace
especially in the general semifinite case that we were interested in.
Specifically for $T\in{\mathcal L}^{(p,\infty)}$ we asked the question of
when the functional
$A\to \tr(AT^s)$ on $\mathcal N$
may be used to calculate the Dixmier trace $\tr_\omega(AT^p)$.
The strongest known result of which we were aware is contained
in Proposition 4 page 306 of \cite{Co4}: for compact operators $T\geq 0$ 
whose singular values $\mu_n(T)$
satisfy $\sum_{n=0}^{N-1}\mu_n(T)=O(\log(N))$, when either
$\lim_{s\to 1}(s-1)\tr(T^s)$ or
$\lim_{N\to\infty} (\log N)^{-1}\sum_{n=0}^{N-1}\mu_n(T)$
exists they both do and are equal. While the somewhat nontrivial proof is
not given there, it does follow as Connes states from the Hardy-Littlewood
Tauberian Theorem (Theorem 98 of \cite{H} is a good reference).
In the PhD thesis of Prinzis \cite{P} an extension of this latter
result was claimed in the type $II$ setting however, the proof was
flawed. Additional interesting information is contained
in \cite{Co4} (page 563) where the Dixmier trace
is expressed in terms of the asymptotics of the trace
of the `heat operator' $e^{\lambda^{-2/p}T^{-2}}$
as $\lambda \to \infty$. Subsequently a proof for this result due to Connes 
was published for $p>1$ in \cite{GVF}.

Our aim in this paper is to prove the strongest possible
theorem relating the zeta function, the asymptotics
of the trace of the heat operator and the Dixmier trace
in both the type $I$ and type $II$ setting
of ${\mathcal L}^{(p,\infty)}$ summable (Breuer-)Fredholm
modules ($1\leq p <\infty$). We obtain the most general results possible
in the most general semifinite case without assuming that any of the above
limits exist. To do this we need a rather novel approach
to the Dixmier trace which we explain in the first section.
The essence of our approach is contained in
Theorem 1.5 where we observe that
there are really two
Dixmier traces, one which might naturally be regarded as
being constructed from an
invariant mean on $L^\infty(\IR)$ (with the additive group
structure on $\IR$) and the other an
invariant mean on $L^\infty(\IR^*_+$) with the multiplicative
group structure on $\IR^*_+$.
The former trace is natural from
the viewpoint of the zeta function
while the latter is that encountered in \cite{Co4}.
Our key observation in Section 3, where we prove the main
Theorems 3.1 and 3.8, is that in order to
calculate the Dixmier trace using the zeta function
 these traces have to be chosen
in pairs related one to the other via the isomorphism
from $\IR$ to $\IR^*_+$ given by the exponential function.

Choose a faithful, normal, semi-finite trace $\tau$ on $\mathcal N$
($\tau$ will be fixed throughout).
Let $D_0$ have resolvent in the ideal of compact operators in $\mathcal N$.
An odd ${\mathcal L}^{(1,\infty)}$
summable unbounded (Breuer)-Fredholm module
for a Banach *-algebra, $\mathcal A$ is a triple $({\mathcal N},{\mathcal A},D_0)$
where ${\mathcal A}\subset{\mathcal N}$ is such that $[a,D_0]$ is
bounded for all $a$ in a dense subalgebra of ${\mathcal A}$ and
$(1+D_0^2)^{-1/2}\in {\mathcal L}^{(1,\infty)}$.
Our main results (in Section 3) concern the asymptotics
of $\tau(A(1+D_0^2)^{-s})$ as $s\to 1/2$ for $A\in \mathcal N$
and how this relates to the Dixmier trace $\tau_\omega(A(1+D^2)^{-1/2})$.
Then in Section 4 we consider the asymptotics of the
trace of the heat semigroup of $D_0^2$ deriving in particular
the formula of (\cite{Co4} p.563) for the Dixmier trace.
Section 5 generalises all of the previous formulae to the case
where $(1+D_0^2)^{-1/2}\in {\mathcal L}^{(p,\infty)}$
with $p>1$.

In Theorem 6.2 we apply our results on the zeta function
approach to the Dixmier trace,
using \cite{CP1} and \cite{CP2}, to derive a general formula
for the Chern character of an ${\mathcal L}^{(1,\infty)}$ summable 
Breuer-Fredholm module $({\mathcal N},{\mathcal A},D_0)$.

In Section 7 we give a brief overview of the  results
in \cite {P} on a Wodzicki residue formula for the Dixmier
trace of pseudo-differential operators tangential
to a minimal ergodic action of $\IR^n$ on a compact space.
Our aim here is to show how the results of the earlier
sections may be used to overcome a technical difficulty
in Prinzis' approach.

Section 8 contains our short proof of the theorem of Lesch
giving the index of a generalised Toeplitz operator
associated with an action of $\IR$ on a $C^*$-algebra
equipped with an invariant trace. The argument depends
in an essential way on our results in Section 3
on the type $II$ Dixmier trace and zeta function
and shows that the index
theory of Toeplitz operators with noncommutative
symbol is a corollary of results in noncommutative geometry.

\subsection{\bf Generalities on singular traces}

We have two groups, the additive group $\IR$ and the multiplicative
group $\IR^*_+$ of positive reals.
The exponential map and the $\log$ are mutually inverse isomorphisms
between these groups.
Notice that $\exp$ takes translation by $a\in \IR$ to dilation by $\exp(a)
\in \IR^*_+$
and dilation by $b\in\IR^*_+$ to the transformation $x\mapsto x^b$
on $\IR^*_+$. Let $G_1$ and $G_2$ be given by
taking the semidirect product of the group $\IR$
and dilations and the semidirect product of the group of
powers with $\IR^*_+$  respectively. That is,
$G_1$ is the set $\IR\times\IR^*_+$ with multiplication:
$$ (a,s)(b,t) = (a+sb,st).$$
While, $G_2$ is the set $\IR^*_+ \times \IR^*_+$ with multiplication:
$$ (s,t)(x,y) = (sx^t,ty).$$
Then, $\exp$ and $\log$ induce mutually inverse isomorphisms
of $G_1$ and $G_2$. For example, the isomorphism $G_1\to G_2$ is given by:
$$ (a,s)\mapsto (exp(a),s):G_1\to G_2.$$

\begin{defn}
We define the isomorphism $L: L^\infty(\IR)\to L^\infty(\IR^*_+)$
by $L(f) = f\circ \log$. We also define the Hardy and Cesaro means (transforms)
on $L^\infty(\IR)$ and $L^\infty(\IR^*_+)$, respectively by:
$$H(f)(u)=\frac{1}{u}\int_0^u f(v) dv\;for\;f\in L^\infty(\IR),\ u\in \IR$$
and,
$$M(g)(t)=\frac{1}{\log t}\int_1^t g(s) \frac{ds}{s}\;for\;g\in
L^\infty(\IR^*_+),\ t>0.$$
We refer to $H$
as the mean for the additive group $\IR$.
\end{defn}

Then a brief calculation yields for $g\in  L^\infty(\IR^*_+)$,
$$LHL^{-1}(g)(r) = \frac{1}{\log r}\int_0^{\log r} g(e^u)du
= \frac{1}{\log r}\int_1^{r} g(v)\frac{dv}{v}= M(g)(r).$$
So indeed $L$ intertwines the two means.

\begin{defn}
We also define the following families of self-maps on these $L^\infty$ spaces:
let $T_b$ denote translation by $b\in \IR$, $D_a$ denote
dilation by $a\in \IR^*_+$ and let $P^a$ denote exponentiation
by $a\in \IR^*_+$. That is,
\begin{eqnarray*}
T_b(f)(x) &=& f(x+b)\;for\; f\in L^\infty(\IR),\\
D_a(f)(x) &=& f(ax)\;for\;f\in L^\infty(\IR),and\\
P^a(f)(x) &=& f(x^a)\;for\;f\in L^\infty(\IR^*_+).
\end{eqnarray*}
\end{defn}

Some of the basic relations between these $L^\infty$ spaces and their
self-maps are provided for easy access by the following proposition.

\begin{prop}
$L^\infty(\IR)$ together with the self-maps, $D_a$, $T_b$, and $H$
($a>0,b\in\IR$) is related to $L^\infty(\IR^*_+)$ together with the
self-maps, $P^a$, $D_a$, and $M$ ($a>0$) via the isomorphism
$$L: L^\infty(\IR)\to L^\infty(\IR^*_+)$$ and the following identities:\\
(1) $LD_aL^{-1} = P^a$ for $a>0$,\\
(2) $LT_bL^{-1} = D_{exp(b)}$ for $b\in\IR$ (or $LT_{log(a)}L^{-1} = D_a$ for
$a>0$),\\
(3) $LHL^{-1} = M$,\\
(4) $D_aH = HD_a$ and $P^aM = MP^a$ for $a>0,$\\
(5) $\lim_{t\to\infty} (HT_b-T_bH)f(t)=0$ for $f\in L^\infty(\IR)$ and
    $b\in\IR,$\\
(6) $\lim_{t\to\infty} (MD_a-D_aM)f(t)=0$ for $f\in L^\infty(\IR^*_+)$ and
    $a>0.$
\end{prop}

\begin{proof}
We have already shown (3). 
The calculations for (1), (2), and (4) are equally
straightforward. To see (5), take $b\in\IR$ and $f\in L^\infty(\IR)$, then:
\begin{eqnarray*}(HT_b-T_bH)f(t)&=&\frac{1}{t}\int_0^t f(x+b)dx -
\frac{1}{t+b}\int_0^{t+b}f(x)dx\\
&=&\frac{1}{t}\int_b^{t+b} f(x)dx - \frac{1}{t+b}\int_0^{t+b}f(x)dx\\
&=&\left(\frac{1}{t}-\frac{1}{t+b}\right)\int_b^{t+b} f(x)dx -
\frac{1}{t+b}\int_0^bf(x)dx\\
&=&\frac{b}{t(t+b)}\int_b^{t+b} f(x)dx - \frac{1}{t+b}\int_0^bf(x)dx.
\end{eqnarray*}

In absolute value this is less than or equal to:
$$\frac{||f||\cdot|b|}{|t+b|} + \frac{||f||\cdot|b|}{|t+b|} =
\frac{2||f||\cdot|b|}{|t+b|},$$
which vanishes as $t\to\infty.$

The proof of (6) is similar.
\end{proof}

We give $G_1$ and $G_2$ the discrete topology to simplify
the discussion and note that they are amenable being extensions of
one abelian group by  a second.
They act as groups of homeomorphisms of $\IR$ and $\IR^*_+$
respectively via $\tilde\alpha_{a,s}(y)=a+sy$ for $(a,s)\in G_1,\;y\in\IR$ and
$\alpha_{s,t}(x)=sx^t$ for $(s,t)\in G_2,\;x\in\IR^*_+$.
Furthermore there are actions of the groups $G_1$ and $G_2$
on $L^\infty(\IR)$ and $L^\infty(\IR^*_+)$.
These actions are generated by $\{T_b,D_a\ \vert\ b\in\IR, a\in \IR^*_+\}$
 in the case of $G_1$
and $\{D_a, P^c\ \vert \ a,c\in \IR^*_+\}$
 in the case of $G_2$ and $L$ intertwines these
actions.
Thus we have actions
$$G_1 \times L^\infty(\IR)^*\to L^\infty(\IR)^*\;
and\;G_2 \times L^\infty(\IR^*_+)^*\to L^\infty(\IR^*_+)^* $$
given respectively by
$$[(a,s), \tilde\omega]\mapsto\tilde\alpha_{a,s}^*(\tilde\omega)\;where\;
\tilde\alpha_{a,s}^*(\tilde\omega)(f)=
\tilde\omega(f\circ\tilde\alpha_{a,s}^{-1})
\;for\; (a,s)\in G_1, \tilde\omega\in L^\infty(\IR)^*\ ,f\in L^\infty(\IR)$$
and
$$[(s,t),\omega]\mapsto \alpha_{s,t}^*(\omega)\;where\;\alpha_{s,t}^*(\omega)(f)
=\omega(f\circ\alpha_{s,t}^{-1})
\;for\; (s,t)\in G_2, \omega\in  L^\infty(\IR^*_+)^*, f\in L^\infty(\IR^*_+).$$

These are weak$^*$-continuous actions because, for example,
if $\omega_\beta\to \omega$ is a net in $L^\infty(\IR^*_+)^*$ then
$$|\alpha_{s,t}^*(\omega_\beta)(f)-\alpha_{s,t}^*(\omega)(f)|
=|\omega_\beta(f\circ \alpha_{s,t}^{-1})-
\omega(f\circ \alpha_{s,t}^{-1})|\to 0$$
as $f\circ \alpha_{s,t}^{-1}\in L^\infty(\IR^*_+)$.

   If we use these remarks together with the previous proposition we obtain the
following.

\begin{prop}
Given any continuous functional $\tilde\omega$
on $L^\infty(\IR)$ which is invariant under $H$ and $G_1$
then  $\tilde \omega\circ L^{-1}$ is
a continuous functional on  $L^\infty(\IR^*_+)$ invariant under $M$ and $G_2$.
Conversely, composition with $L$ converts an $M$
and $G_2$ invariant continuous functional on $L^\infty(\IR^*_+)$
into an $H$ and $G_1$ invariant continuous functional on  $L^\infty(\IR)$.
\end{prop}

\subsection{\bf Existence of invariant singular traces}

We denote by $C_0(\IR^*_+)$ the continuous functions
on $\IR^*_+$ vanishing at infinity. Our aim in this subsection is to
prove the following result.

\begin{thm}
There exists a state $\omega$ on $L^\infty(\IR^*_+)$
satisfying the following conditions:\\
(1) $\omega(C_0(\IR^*_+)) \equiv 0$.\\
(2) If $f$ is real-valued in $L^\infty(\IR^*_+)$
then
$$ess\ liminf_{t\to\infty} f(t) \leq\omega(f)\leq ess\ limsup_{t\to\infty}
f(t).$$
(3) If the essential support of $f$ is compact then $\omega(f)=0.$\\
(4) For all $c\in \IR^*_+$, $\omega(D_cf)=\omega(f)$ for all
$f\in L^\infty(\IR^*_+)$.\\
(5) For all $a \in\IR^*_+$ and all $f\in L^\infty(\IR^*_+)$
 $\omega(P^af)=\omega(f)$.\\
(6) For all  $f\in L^\infty(\IR^*_+)$, $\omega(Mf)=  \omega(f)$.

\end{thm}

Using the preceding proposition we obtain the following:

\begin{cor}
There exists a state $\tilde\omega$ on $L^\infty(\IR)$
satisfying the following conditions:\\
(1) $\tilde\omega(C_0 (\IR)) \equiv 0$.\\
(2) If $f$ is real-valued in $L^\infty(\IR)$
then
$$ess\ liminf_{t\to\infty} f(t) \leq\tilde\omega(f)\leq ess\ limsup_{t\to\infty}
f(t).$$
(3) If the essential support of $f$ is compact then $\tilde\omega(f)=0.$\\
(4) For all $c\in \IR$, $\tilde\omega(T_cf)=\tilde\omega(f)$ for all
$f\in L^\infty(\IR)$.\\
(5) For all $a \in\IR^*_+$ and all $f\in L^\infty(\IR)$
 $\tilde\omega(D_af)=\tilde\omega(f)$.\\
(6) For all  $f\in L^\infty(\IR)$, $\tilde\omega(Hf)= \tilde\omega(f)$.
\end{cor}

Notice that $L$ sends $C_0 (\IR))$ into 
 $C_0 (\IR_+^*))$.
Also,
 we observe that condition 2 of the corollary is equivalent to the statement
that if $f\in L^\infty(\IR)$ is continuous and
$\lim_{|t|\to\infty}f(t)$ exists then 
$\tilde\omega(f)= \lim_{|t|\to\infty}f(t)$.
The rest of this subsection will be devoted to the proof of the
theorem.
Introduce the set $S$ consisting of all positive functionals
$\omega\in L^\infty(\IR^*_+)^* $
normalised so that $\omega(1)=1$
and such that condition 1 of the theorem holds.

Clearly $S$ is a convex and weak$^*$ closed subset of the unit ball.
Moreover $S$ is non-empty as we can define
$\omega\in (C[0,\infty])^*$ by $\omega(f)=f(\infty)$
then $\omega$ is positive and $\omega(1)=||\omega||=1$.
So extending $\omega$ to $\tilde \omega \in L^\infty(\IR^*_+)^*$
by Hahn-Banach yields a non-trivial element of $S$
(note that positivity of the extension is well known,
for example see Theorem 4.3.2 of \cite{KR}).

It is straightforward to verify $G_2$ acts affinely (i.e. preserving
convex combinations) on $S$
by restriction of the dual action on  $L^\infty(\IR^*_+)^*$.
As we have remarked earlier the action is weak$^*$ continuous and
$G_2$ is amenable since it is the extension of an abelian group by an
abelian group (and so too is $G_1$). Hence
by Rickert's Theorem \cite{G} there is a fixed point $\omega_0$
for this action.
This fixed point satisfies conditions (1), (4) and (5) of the theorem.
Condition (3) holds because if $f\geq 0$ and has compact support then
there is a continuous function $g\geq f\;a.e.$ with $g(\infty)=0$ and so,
$$ 0\leq\omega_0(f)\leq\omega_0(g)=0.$$
To see that $\omega_0$ satisfies condition (2),
let $f$ be real-valued and let $C$
denote the $ess\limsup_{t\to\infty}f(t)$.
Then for each $\epsilon>0$ there exists a function $g$ with
the support of $g$ compact and $(f-g)\leq C +\epsilon$ a.e.
Then $\omega_0(f)\leq C +\epsilon$. Similarly
$\omega_0(f)$ is bounded below by the  essential $\liminf_{t\to\infty}f(t)$

Let $M^*$ denote the linear map on  $L^\infty(\IR^*_+)^*$
given by $M^*\omega(f)=\omega(Mf)$.
Finally, to prove (6) we note first that $M$ leaves $C_0(\IR^*_+)$
and the constant function invariant and hence that $M^*$ leaves $S$ invariant.
By Proposition 1.3 (part (6) and the second half of part (4)), we see that
the action of $M^*$ ({\bf on} $S$!)
commutes with the dual actions of the generators,
$D_a$ and $P^a$ of $G_2$ ({\bf on} $S$!).
It follows then that for any fixed point $\omega_0$ of the $G_2$ action,
$\omega_0\circ M^*$ is another fixed point of the $G_2$ action on $S$.
In other words, $M^*$ leaves the set of $G_2$ fixed points of
$S$ invariant. Thus $M^*$ leaves the set of functionals in $S$
satisfying conditions (1) to (5) invariant.
 The collection of fixed points for
$G_2$ is clearly a weak-* compact convex set invariant under the (affine)
action of $M^*$. It follows from the Kakutani-Markov Theorem \cite{E}
that $M^*$ itself has a fixed
point in this subset which is therefore a functional satisfying conditions
(1) to (6) of the theorem completing the proof.

\noindent{\bf Remarks}:
The spirit of the approach of this section goes back to Dixmier
\cite{Dix1}. The approach of Connes \cite{Co4} is different in
a slightly subtle way which we will not go into fully here.
Suffice to say that \cite{Dix1} uses dilation invariant functionals
from the start while \cite{Co4} uses the Cesaro mean to obtain
a dilation invariant functional (that is, starting from a state
$\omega$ on  $ L^\infty(\IR^*_+)^* $ one observes that $\omega\circ M^*$
is dilation invariant). This difference is important to us in Sections 5
and 6.

\subsection{\bf Notation}
We are interested in  certain ideals of operators in
the von Neumann algebra $\mathcal N$  defined using our faithful, normal,
semifinite trace $\tau$.

\begin{defn}
 If $S\in\mathcal N$ the {\bf t-th generalized
singular value} of S for each real $t>0$ is given by
$$\mu_t(S)=\inf\{||SE||\ \vert \ E \text{ is a projection in }
{\mathcal N} \text { with } \tau(1-E)\leq t\}.$$
\end{defn}

We will mostly explain
the results we need about these singular values
later in the text although a full exposition is
contained in \cite{F} and 
\cite{FK}. We write
$T_1\prec\prec T_2$ to mean that $\int_0^t \mu_s(T_1)ds\leq
\int_0^t \mu_s(T_2)ds$ for all $t>0$. 

\begin{defn}
If $\mathcal I$ is a $*$-ideal in $\mathcal N$
which is complete in a norm $||\cdot||_{\mathcal I}$
then we will call $\mathcal I$ an {\bf invariant
operator ideal} if\\
(1) $||S||_{\mathcal I}\geq ||S||$ for all $S\in \mathcal I$,\\
(2) $||S^*||_{\mathcal I} = ||S||_{\mathcal I}$ for all $S\in \mathcal I$,\\
(3) $||ASB||_{\mathcal I}\leq ||A|| \:||S||_{\mathcal I}||B||$ for all
$S\in \mathcal I$, $A,B\in \mathcal N$.\\
Since $\mathcal I$ is an ideal in a von Neumann algebra, it follows
from I.1.6, Proposition 10 of \cite{Dix} that if $0\leq S \leq T$
and $T \in {\mathcal I}$, then $S \in {\mathcal I}$ and $||S||_{\mathcal I}
\leq ||T||_{\mathcal I}$. Much more is true, especially in the type $I$ case
but we shall not need it here, see \cite{GK}.
\end{defn}

The main examples of such ideals that we consider in this paper
are the spaces
$${\mathcal L}^{(1,\infty)}({\mathcal N})=\left\{T\in{\mathcal N}\ |\ \Vert T\Vert_{_{{\mathcal L}^{(1,\infty)}}}:=   \sup_{t> 0}
\frac{1}{\log(1+t)}\int_0^t\mu_s(T)ds<\infty\right\}.$$
and with $p>1$,
$$\psi_p(t)=\left\{\begin{array}{ll} t & \mbox{for } 0\leq t\leq 1\\
                                     t^{1-\frac{1}{p}} & \mbox{for } 1\leq t
\end{array}\right.$$
$${\mathcal L}^{(p,\infty)}({\mathcal N})=\left\{T\in{\mathcal N}\ |\ \Vert T\Vert_{_{{\mathcal L}^{(p,\infty)}}}:= \sup_{t> 0}
\frac{1}{\psi_p(t)}\int_0^t\mu_s(T)ds<\infty\right\}.$$

There is also the equivalent definition
$${\mathcal L}^{(p,\infty)}({\mathcal N})=\left\{T\in{\mathcal N}\ |\ \sup_{t> 0}
\frac{t}{\psi_p(t)}\mu_t(T)<\infty\right\}.$$

It is well-known (see e.g. \cite{GK},\cite{Co4}) that for $T_1\in {\mathcal N}$, $T_2\in {\mathcal L}^{(p,\infty)}({\mathcal N})$, $p\in [1,\infty)$, 
the condition $T_1\prec\prec T_2$ implies that $T_1\in {\mathcal L}^{(p,\infty)}({\mathcal N})$.

As we will not change ${\mathcal N}$ throughout the paper we will suppress
the $({\mathcal N})$ to lighten the notation. On this point however
the reader should note
that ${\mathcal L}^{(p,\infty)}$ is often taken to mean an ideal in
the algebra $\widetilde{\mathcal N}$ of measurable operators affiliated
to ${\mathcal N}$. Our notation is however consistent with that of \cite{Co4}
in the special case ${\mathcal N}={\mathcal B}({\mathcal H})$.

For most of the paper
$T$ is a positive operator in ${\mathcal L}^{(1,\infty)}$.
There is a map from  the positive operators in ${\mathcal L}^{(1,\infty)}$
to $L^\infty[0,\infty)$ given by
$T\to f_T$ where $f_T(t)=\frac{1}{\log(1+t)}\int_0^t \mu_s(T) ds$.
We may extend $f_T$ to all of $\IR$ by defining it to be
zero on the negative reals. Depending on the circumstances
we can thus regard $f_T$ as either an element of
$L^\infty(\IR)$ or $L^\infty(\IR^*_+)$.

Henceforth we use the notation $\tau_\omega(T)$
for $\omega(f_T)$ where
$f_T(t)=\frac{1}{\log(1+t)}\int_0^t \mu_s(T) ds$
and $\omega\in L^\infty(\IR^*_+)^*$
satisfies the conditions of Theorem 1.5.
We also write
$$\tau_\omega(T)=
\omega-\lim_{t\to\infty}\frac{1}{\log(1+t)}\int_0^t \mu_s(T) ds.$$
It follows from \cite{Co4}, IV.2.$\beta$ (see also \cite{DPSS}, Example 2.5) 
that $\tau_\omega(\cdot)$ is additive and positively homogeneous on the 
positive part of ${\mathcal L}^{(1,\infty)}$ and hence extends to a positive 
linear functional  on ${\mathcal L}^{(1,\infty)}$ (again denoted by $\tau_\omega$).
It is in fact an example of a singular trace on
$\mathcal N$ ({\it cf} the discussion in \cite{Co4} and \cite {DPSS})

\section{\bf Preliminary results}

It is useful to have an estimate on the singular values
of the operators in ${\mathcal L}^{(1,\infty)}$.

\begin{lemma}
For $T\in{\mathcal L}^{(1,\infty)}$ positive
there is a constant $K>0$ such that for each
$p\geq 1$,
$$\int_0^t\mu_s(T)^pds \leq K^p\int_0^t\frac{1}{(s+1)^p}ds.$$
\end{lemma}

\begin{proof}  By \cite{FK}, 
Lemma 2.5 (iv),  for all $0\leq T\in {\mathcal N}$ and 
all continuous increasing functions $f$ on $[0,\infty)$ with $f(0)\ge 0$, 
we have $\mu_s(f(T))=f(\mu_s(T))$ for all $s>0$. 
Combining this fact with well-known result of Hardy-Littlewood-P\' olya 
(see e.g. \cite{F}, Lemma 4.1), we see that $T_1\prec\prec T_2$, 
$0\leq T_1,T_2\in {\mathcal N}$ implies $T_1^p\prec\prec T_2^p$ for 
all $p\in (1,\infty)$. Now, by definition of ${\mathcal L}^{(1,\infty)}$ 
the singular values of $T$ satisfy $\int_0^t\mu_s(T)ds=O(\log t)$ so that
for some $K>0$,
$$\int_0^t\mu_s(T) ds \leq K \int_0^t\frac{1}{(s+1)}ds,\quad \forall t>0.$$
In other words $\mu_s(T)\prec\prec K/(1+s)$ and the assertion of lemma 
follows immediately.
\end{proof}

\begin{thm}(weak$^*$-Karamata theorem)
Let $\tilde\omega\in L^\infty(\IR)^*$
be a dilation invariant state and
let $\beta$ be a real valued, increasing,
right continuous function on $\IR_+$
which is zero at zero
and such that the integral
$h(r)=\int_0^\infty e^{-\frac{t}{r}}d\beta(t)$
converges for all $r>0$ and
$C=\tilde\omega-\lim_{r\to\infty} \frac{1}{r}h(r)$ exists.
Then
$$\tilde\omega-\lim_{r\to\infty} \frac{1}{r}h(r)
=\tilde\omega-\lim_{t\to\infty}\frac{\beta(t)}{t}.$$
\end{thm}

\noindent{\bf Remark}:
The classical Karamata theorem states, in the notation of the
theorem, that if the ordinary limit
$\lim_{r\to\infty} \frac{1}{r}h(r)=C$ exists
then $C=\lim_{t\to\infty}\frac{\beta(t)}{t}$.
The proof of this classical result is obtained by replacing,
in the proof of Theorem 2.2,
$\tilde\omega-\lim$ throughout by the ordinary limit.

\begin{proof}
Let
$$g(x)=\left\{\begin{array}{ll} x^{-1} & \mbox{for } e^{-1}\leq x\leq 1\\
                                0       & \mbox{for } 0\leq x< e^{-1}
\end{array}\right.$$
so that $g$ is right continuous at $e^{-1}$.
Then for $r>0$, $t\to e^{-t/r}g(e^{-t/r})$ is left continuous
at $t=r$. Thus the Riemann-Stieltjes integral
$\int_0^\infty e^{-t/r}g(e^{-t/r})d\beta(t)$
exists for each $r>0$.
We claim that for any polynomial $p$
$$\tilde\omega-\lim_{r\to\infty}\frac{1}{r}\int_0^\infty
e^{-t/r}p(e^{-t/r})d\beta(t) =C\int_0^\infty
e^{-t}p(e^{-t})dt.$$ To see this first compute
for $p(x)=x^n$,
$$\frac{1}{r}\int_0^\infty e^{-t/r}e^{-nt/r}d\beta(t)
=\frac{1}{r}\int_0^\infty e^{-(n+1)t/r}d\beta(t).$$
Therefore
$$\frac{1}{n+1}\tilde\omega-\lim_{r\to\infty}\frac{1}{r/(n+1)}
\int_0^\infty e^{-(n+1)t/r}d\beta(t) =\frac{C}{n+1}$$
by dilation invariance of $\tilde\omega$.
Thus $$\tilde\omega-\lim_{r\to\infty}
\frac{1}{r}\int_0^\infty e^{-t/r}e^{-nt/r}d\beta(t)
=C\int_0^\infty
e^{-t}(e^{-t})^ndt.$$ Since $\tilde\omega$ is linear
the claim follows for all $p$.

Choose sequences of polynomials
$\{p_n\},\ \{P_n\}$ such that for
all $x\in [0,1]$
$$-1\leq p_n(x)\leq g(x) \leq P_n(x)\leq 3$$
and such that $p_n$ and $P_n$ converge a.e. to $g(x)$.
Then since $\tilde\omega$ is positive it preserves order:
$$C\int_0^\infty
e^{-t}p_n(e^{-t})dt=\tilde\omega-\lim_{r\to\infty}\frac{1}{r}\int_0^\infty
e^{-t/r}p_n(e^{-t/r})d\beta(t)\leq \tilde\omega-\lim_{r\to\infty}\frac{1}{r}\int_0^\infty
e^{-t/r}g(e^{-t/r})d\beta(t)$$
$$\leq\ldots\leq C\int_0^\infty
e^{-t}P_n(e^{-t})dt.$$
By the Lebesgue Dominated Covergence Theorem both
$\int_0^\infty
e^{-t}p_n(e^{-t})dt$ and $\int_0^\infty
e^{-t}P_n(e^{-t})dt$ converge to $\int_0^\infty
e^{-t}g(e^{-t})dt$ as $n\to\infty$.
But a direct calculation yields $\int_0^\infty
e^{-t}g(e^{-t})dt=1$ and
$$\int_0^\infty e^{-t/r}g(e^{-t/r})d\beta(t)=\beta(r).$$
Hence
$$C=
\tilde\omega-\lim_{r\to\infty}\frac{1}{r}\int_0^\infty
e^{-t/r}g(e^{-t/r})d\beta(t)
=\tilde\omega-\lim_{r\to\infty}\frac{\beta(r)}{r}.$$
\end{proof}

\noindent Recall that for any $\tau$-measurable operator $T$, the distribution
function of $T$ is defined by
$$
\lambda _t(T):=\tau (\chi_{(t,\infty)}(|T|)),\quad t>0,
$$
where $\chi_{(t,\infty)}(|T|)$ is the spectral projection of $|T|$
corresponding to the interval $(t,\infty)$ (see [FK]). By Proposition 2.2
of [FK],
$$\mu_s(T)=\inf\{t\ge 0\ :\ \lambda_t(T)\leq s\}
$$
we infer that for any  $\tau$-measurable operator $T$, the distribution
function $\lambda _{(\cdot)}(T)$ coincides with the (classical) distribution
function of $\mu_{(\cdot)}(T)$. From this formula and the fact that $\lambda$
is right-continuous, we can easily see that for $t>0$, $s>0$
$$s\geq\lambda_t\Longleftrightarrow \mu_s\leq t.$$
Or equivalently,
$$s < \lambda_t \Longleftrightarrow \mu_s > t.$$
Using Remark 3.3 of [FK] this implies that:
$$ \int_0^{\lambda _t} \mu_s(T)ds=\int_{[0,\lambda _t)} \mu_s(T)ds=
\tau (|T|\chi_{(t,\infty)}(|T|)),\quad t>0.\eqno (*)$$

\begin{lemma}
For $T\in{\mathcal L}^{(1,\infty)}$ and
$C>\Vert T\Vert _{{\mathcal L}^{(1,\infty)}}$ we have eventually
$$
\lambda _{\frac{1}{t}}(T)\leq Ct\log t.
$$
\end{lemma}

\begin{proof} Suppose not and there exists $t_n\uparrow \infty$ such that
$\lambda _{\frac{1}{t_n}}(T)> Ct_n\log t_n$ and so for $s\leq Ct_n\log t_n$
we have $\mu_s(T)\geq\mu_{Ct_n\log t_n}(T) > \frac{1}{t_n}$. Then for sufficiently
large $n$
$$
\int_0^{ Ct_n\log t_n}\mu_s(T)ds> \frac{1}{t_n}\cdot  Ct_n\log t_n=
C\log t_n.
$$
Choose $\delta>0$ with $C-\delta>\Vert T\Vert_{{\mathcal L}^{(1,\infty)}}$.
 Then for sufficiently
large $n$
$$
C\log t_n=(C-\delta)\log t_n +\delta\log t_n>
\Vert T\Vert _{{\mathcal L}^{(1,\infty)}}\log( Ct_n)
+\Vert T\Vert _{{\mathcal L}^{(1,\infty)}}\log(\log (t_n+1))
$$
$$
=\Vert T\Vert _{{\mathcal L}^{(1,\infty)}}\log( Ct_n\log (t_n+1)).$$
This is a contradiction with the inequality
$\int_0^{t}\mu_s(T)ds\leq \Vert T\Vert_{{\mathcal L}^{(1,\infty)}}\log (t+1)$,
which holds for any $t>0$ due to the
definition of the norm in ${\mathcal L}^{(1,\infty)}$.
\end{proof}

\noindent An assertion somewhat similar to Proposition 2.4
 below was formulated
in [P] and supplied with an incorrect proof. We use a different approach.
\bigskip

\begin{prop}
For $T\in{\mathcal L}^{(1,\infty)}$ positive let $\omega$ be a $G_2$ invariant
state on  $L^\infty(\IR^*_+)$. For every $C>0$
$$
\tau_\omega(T)=
\omega-\lim_{t\to\infty}\frac{1}{\log(1+t)}\int_0^t \mu_s(T) ds
=\omega-\lim_{t\to\infty}\frac{1}{\log(1+t)}
\tau(T\chi_{(\frac{1}{t},\infty)}(T))
$$
$$
=\omega-\lim_{t\to\infty}\frac{1}{\log(1+t)}\int_0^{Ct\log t}\mu_s(T) ds
$$
and if one of the $\omega-$limits is a true limit then so are the others.

\end{prop}

\begin{proof} We first note that
$$
\int_0^t\mu_s(T)ds\leq \int_0^{\lambda _{\frac{1}{t}}(T)} \mu_s(T)ds+1,
\quad t>0.
$$
Indeed, the inequality above holds trivially if
$t\leq\lambda _{\frac{1}{t}}(T)$. If $t> \lambda _{\frac{1}{t}}(T)$, then
$$
\int_0^t\mu_s(T)ds=\int_0^{\lambda _{\frac{1}{t}}(T)} \mu_s(T)ds+\int_
{\lambda _{\frac{1}{t}}(T)}^t \mu_s(T)ds.
$$
Now $s> \lambda _{\frac{1}{t}}(T)$ implies that $\mu_s(T)\leq \frac{1}{t}$
so we have
$$
\int_0^t\mu_s(T)ds\leq \int_0^{\lambda _{\frac{1}{t}}(T)} \mu_s(T)ds+\frac{1}{t}(t-\lambda _{\frac{1}{t}}(T))\leq \int_0^{\lambda _{\frac{1}{t}}(T)} \mu_s(T)ds+1.
$$
Using this observation and lemma  above we see that for
$C>\Vert T\Vert _{{\mathcal L}^{(1,\infty)}}$ and any fixed
$\alpha >1$  eventually
$$
\int_0^t\mu_s(T)ds \leq \int_0^{\lambda _{\frac{1}{t}}(T)} \mu_s(T)ds+1\leq \int_0^{Ct\log t} \mu_s(T) ds +1 \leq\int_0^{t^\alpha}
\mu_s(T) ds+1
$$
and so eventually
$$
\frac{1}{\log(1+t)}\int_0^t\mu_s(T)ds
\leq\frac{1}{\log(1+t)}( \int_0^{\lambda _{\frac{1}{t}}(T)} \mu_s(T)ds+1)
\leq\frac{1}{\log(1+t)}(\int_0^{Ct\log t} \mu_s(T) ds +1)
$$
$$
\leq\frac{\log(1+t^\alpha)}{\log(1+t)\log(1+t^\alpha)}
(\int_0^{t^\alpha} \mu_s(T) ds+1).
$$
Taking the $\omega$-limit we get
$$
\tau_\omega(T) \leq\omega-\lim_{t\to\infty}\frac{1}{\log(1+t)}
\int_0^{\lambda _{\frac{1}{t}}(T)} \mu_s(T)ds
\leq\omega-\lim_{t\to\infty}\frac{1}{\log(1+t)}\int_0^{Ct\log t} \mu_s(T) ds
$$
$$
\leq\omega-\lim_{t\to\infty}\frac{\alpha}{\log(1+t^\alpha)}
\int_0^{t^\alpha} \mu_s(T) ds=\alpha\tau_\omega(T)
$$
where the last line uses $G_2$ invariance.
Since this holds for all $\alpha>1$ and using $(*)$ we get the conclusion
for $\omega$-limits and $C>\Vert T\Vert _{{\mathcal L}^{(1,\infty)}}$.
The assertion for an arbitrary $0<C\leq \Vert T\Vert _{{\mathcal L}^{(1,\infty)}}$ follows immediately by noting that for $C'> \Vert T\Vert _{{\mathcal L}^{(1,\infty)}}$ one has eventually
$$
\int_0^t\mu_s(T)ds\leq\int_0^{Ct\log t} \mu_s(T) ds\leq\int_0^{C't\log t}
\mu_s(T) ds.
$$

To see the last assertion of the Proposition suppose that
$\lim_{t\to\infty}\frac{1}{\log(1+t)}\int_0^t\mu_s(T)ds=A$
then by the above argument we get
$$
A\leq\liminf_{t\to\infty}\frac{1}{\log(1+t)}
\tau(T\chi_{(\frac{1}{t},\infty)}(T))
\leq\limsup_{t\to\infty}\frac{1}{\log(1+t)}
\tau(T\chi_{(\frac{1}{t},\infty)}(T))
\leq\alpha A
$$
for all $\alpha>1$ and hence
$\lim_{t\to\infty}\frac{1}{\log(1+t)}
\tau(T\chi_{(\frac{1}{t},\infty)}(T))=A$.
On the other hand if the limit
$\lim_{t\to\infty}\frac{1}{\log(1+t)}\tau(T\chi_{(\frac{1}{t},\infty)}(T))$
exists and equals $B$ say then
$$
\limsup_{t\to\infty}\frac{1}{\log(1+t)}\int_0^t \mu_s(T) ds\leq B\leq
\alpha\liminf_{t\to\infty}\frac{1}{\log(1+t)}\int_0^t \mu_s(T) ds
$$
for all $\alpha>1$ and so
$$
\lim_{t\to\infty}\frac{1}{\log(1+t)}\int_0^t \mu_s(T) ds=B
$$
as well. The remaining claims follow similarly.
\end{proof}

\section{\bf The zeta function and the Dixmier trace}

The zeta function of positive $T\in {\mathcal L}^{(1,\infty)}$
is given by
$$\zeta(s)=\tau(T^s)$$
and for $A\in\mathcal N$ we set
$$\zeta_A(s)=\tau(AT^s).$$

We are interested in the asymptotic behaviour of $\zeta(s)$
and $\zeta_A(s)$ as $s\to 1$.

Now it is elementary to see that the discussion of singular
traces is relevant because
by Lemma 2.1 we have for some $K>0$ and all $s>1$
 $$\tau(T^s) =\int_0^\infty \mu_r(T^s) dr=\int_0^\infty \mu_r(T)^s dr
$$
$$\leq \int_0^\infty \frac{K^s}{(1+r)^s} dr = \frac{K^s}{s-1}.$$
{}From this it follows that $\{(s-1)\tau(T^s)\vert \ s>1\}$
is bounded. Now for $A$ bounded
$|(s-1)\tau(AT^s)|\leq ||A||(s-1)\tau(T^s)$
so that $(s-1)\tau(AT^s)$ is also bounded
and hence for any
$\tilde\omega \in L^\infty(\IR)^*$ satisfying conditions (1), (2) and (3)
of Corollary 1.6
$$\tilde\omega-\lim_{r\to \infty}\frac{1}{r}\tau(AT^{1+\frac{1}{r}})
\eqno (3.1)$$
 exists.

Here we think of $r\to \frac{1}{r}\tau(AT^{1+\frac{1}{r}})$
as a function on all of $\IR$ by extending it to be identically zero
for $r<1$. For notational convenience one might like to think of  (3.1)
as  $\tilde\omega-\lim_{s\to 1}(s-1)\tau(AT^s)$ but this of course does
not (strictly speaking) make sense whereas if $\lim_{s\to
  1}(s-1)\tau(AT^s)$
exists then it is $\lim_{r\to\infty}\frac{1}{r}\tau(AT^{1+\frac{1}{r}})$.

In the following theorem we will take $T\in{\mathcal L}^{(1,\infty)}$ positive,
$||T||\leq 1$ with spectral resolution $T=\int \lambda dE(\lambda)$. We would
like to integrate with respect to $d\tau(E(\lambda))$; unfortunately, these
scalars $\tau(E(\lambda))$ are, in general,
 all infinite. To remedy this situation, we instead
must integrate with respect to the increasing (negative) real-valued function
$N_T(\lambda)=\tau(E(\lambda)-1)$ for $\lambda >0$. Away from $0$, the
increments $\tau(\triangle E(\lambda))$ and $\triangle N_T(\lambda)$ are, of
course, identical.

In a recent email, Alain Connes has sent us a proof of the 
more difficult 
implication of Proposition 4 on page 306 of \cite{Co4}. This is the essential 
point in the proof of the second statement of the theorem below for 
$\mathcal N = \mathcal B(\mathcal H)$. While his 
argument is admittedly simpler, it is similar in spirit to the proof below as 
it uses Karamata's approach to the classical Hardy-Littlewood Tauberian
Theorem (Theorem 98 in \cite{H}), as suggested by Connes in \cite{Co4}. 
.

\begin{thm}
For $T\in{\mathcal L}^{(1,\infty)}$ positive, $||T||\leq 1$ and
$\tilde\omega\in L^\infty(\IR)^*$ satisfying all the conditions
of Corollary 1.6, let $\tilde\omega=\omega\circ L$
where $L$ is given in subsection 1.1, then
we have:
$$\tau_\omega(T)=\tilde\omega-\lim\frac{1}{r}\tau(T^{1+\frac{1}{r}}).$$
If $\lim_{r\to\infty}\frac{1}{r}\tau(T^{1+\frac{1}{r}})$ exists then
$$\tau_\omega(T)=\lim_{r\to\infty}\frac{1}{r}\tau(T^{1+\frac{1}{r}})$$
for an arbitrary dilation invariant functional
$\omega\in L^\infty(\IR^*_+)^*$.
\end{thm}

\begin{proof}
By (3.1) we can apply the weak$^*$-Karamata
theorem to $\frac{1}{r}\tau(T^{1+\frac{1}{r}})$.
First write
$\tau(T^{1+\frac{1}{r}})=\int_{0^+}^1 \lambda^{1+\frac{1}{r}} dN_T(\lambda)$.
Thus setting $\lambda=e^{-u}$
$$\tau(T^{1+\frac{1}{r}})=\int_0^\infty e^{-\frac{u}{r}}d\beta(u)$$
where $\beta(u)=\int_u^0 e^{-v}dN_T(e^{-v})=-\int_0^u e^{-v}dN_T(e^{-v})$.
Since the change of variable $\lambda=e^{-u}$ is strictly decreasing, $\beta$
is, in fact, nonnegative and increasing.
By the weak$^*$-Karamata theorem 
 applied to $\tilde\omega\in L^\infty(\IR)^*$
$$\tilde\omega-\lim_{r\to \infty} \frac{1}{r}\tau(T^{1+\frac{1}{r}})
=\tilde\omega-\lim_{u\to\infty}
\frac{\beta(u)}{u}.$$

Next with the substitution $\rho=e^{-v}$ we get:
$$\tilde\omega-\lim_{u\to\infty}\frac{\beta(u)}{u}=
\tilde\omega-\lim_{u\to\infty}
\frac{1}{u}\int^1_{e^{-u}}
\rho dN_T(\rho).  \eqno (3.2)$$

Set $f(u)=\frac{\beta(u)}{u}$.
We want to make the change of variable $u=\log t$
or in other words to consider $f\circ \log = Lf$.
We use the discussion in subsection 1.1
which tells us that if we start with a $G_2$ and $M$ invariant
functional $\omega\in L^\infty(\IR_+^*)^*$
then the functional $\tilde\omega=\omega\circ L$
is $G_1$ and $H$ invariant as required by the
theorem. Then we have
$$\tilde\omega-\lim_{r\to \infty} \frac{1}{r}\tau(T^{1+\frac{1}{r}})
=\tilde\omega-\lim_{u\to\infty}\frac{\beta(u)}{u}=\tilde\omega-\lim_{u\to\infty}
f(u)=\omega-\lim_{t\to\infty}Lf(t)=\omega-\lim_{t\to\infty}
\frac{1}{\log t}\int_{1/t}^1\lambda dN_T(\lambda).$$

Now, by Proposition 2.4
$$\omega-\lim_{t\to\infty}
\frac{1}{\log t}\int_{1/t}^1\lambda dN_T(\lambda)
=\omega-\lim_{t\to\infty} \frac{1}{\log t}\tau(\chi_{(\frac{1}{t},1]}(T)T)
=\tau_\omega(T).$$
This completes the proof of the first part of the theorem.

The proof of the second part is similar. Using the classical
Karamata theorem (see the remark following the statement of Theorem 2.2)
we obtain the following analogue of (3.2):
$$\lim_{r\to\infty}\frac{1}{r}\tau(T^{1+r})=\lim\frac{\beta(u)}{u}=
\lim_{u\to\infty}\frac{1}{u}\int_{e^{-u}}^1\rho dN_T(\rho).$$
Making the substitution $u=\log t$ on the right hand side we have
$$\lim_{u\to\infty}\frac{1}{u}\int_{e^{-u}}^1\rho dN_T(\rho)=
\lim_{t\to\infty}\frac{1}{\log t}\int_{\frac{1}{t}}^1\lambda dN_T(\lambda)
=\tau_\omega(T)$$ where in the last equality we need only
dilation invariance of the state $\omega\in L^\infty(\IR^*_+)^*$
and not the full list of conditions of Corollary 1.6.
\end{proof}

The map on positive $T\in {\mathcal L}^{(1,\infty)}$ to $\IR$
given by $T\to \tau_\omega(T)$ can be extended
by linearity to a $\IC$ valued functional on
all of ${\mathcal L}^{(1,\infty)}$.
Then the functional
$$A\mapsto\tau_\omega(AT) \eqno (**)$$ for $A\in\mathcal N$ and fixed
$T\in{\mathcal L}^{(1,\infty)}$ is well defined.
We intend to study the properties of (**).
Part of the interest in this functional
stems from the following result
as well as the use of the Dixmier trace
in noncommutative geometry \cite{Co4}.

\begin{lemma}Let $T\in{\mathcal L}^{(1,\infty)}$, then\\
(i) For $A\in\mathcal N$  we have
$$\tau_\omega(AT)=\tau_\omega(TA).$$
(ii) Assume that $D_0$ is an unbounded self adjoint
operator affiliated with $\mathcal N$ such that\\
$T=(1+D_0^2)^{-1/2}\in{\mathcal L}^{(1,\infty)}$.
If $[A_j,|D_0|]$ is a bounded operator for $A_j\in{\mathcal N},\ j=1,2$ then
$$\tau_\omega(A_1A_2T)=\tau_\omega(A_2A_1T).$$
\end{lemma}

\begin{proof}(i) This is proposition A.2 of \cite{CM}.
The proof is elementary, first show that
$\tau_\omega(UTU^*)=\tau_\omega(T)$
then use linearity to extend to arbitrary
$T\in{\mathcal L}^{(1,\infty)}$.
Replace $T$ by $TU$ then use linearity again.\\
(ii) We remark that  $[A_j,|D_0|]$
defining a bounded operator means that the $A_j$ leave $dom(|D_0|)=dom(D_0)$
invariant and that $[A_j,|D_0|]$ is bounded on this domain (see [BR] 3.2.55
and its proof for equivalent but seemingly weaker conditions).
As $|D_0|-(1+D_0^2)^{1/2}$ is bounded, $[A_j, (1+D_0^2)^{1/2}]$ defines a
bounded operator whenever
$[A_j,|D_0|]$ does. As $T^{-1} = (1+D_0^2)^{1/2}$ and $T:\mathcal H \to
dom(T^{-1})$, we see that the formal calculation:
$$[A_j,T]=A_jT-TA_j=T(T^{-1}A_j-A_jT^{-1})T=T[T^{-1},A_j]T$$
makes sense as an everywhere-defined operator on $\mathcal H$. That is,
$$[A_j,T]=T[(1+D_0^2)^{1/2},A_j]T\in({\mathcal L}^{(1,\infty)})^2\subseteq
{\mathcal L}^1.$$ Then
we have, using part (i),
$$\tau_\omega(A_1A_2T)=\tau_\omega(A_2A_1T)-\tau_\omega([A_1,T]A_2).$$
So then
$$\tau_\omega(A_1A_2T)
=\tau_\omega(A_2A_1T)-\tau_\omega(T[(1+D_0^2)^{1/2},A_1]TA_2).$$
Since the operator in the last term is trace class we are done.
\end{proof}

As a corollary of this lemma we see that
 (**) can be used to define a trace on certain
subalgebras of ${\mathcal N}$.
We aim to give several formulas for it.
The first involves the zeta function.
We begin with some preliminary lemmas.

\begin{lemma}
Let $T\geq 0$, $b\geq 0$ be bounded operators\\
(i) If $||b||\leq M$ then for any $1\leq s <2$
$$(b^{1/2}Tb^{1/2})^s\leq M^{s-1}b^{1/2}T^sb^{1/2}.$$
(ii) If $m>0$, $\bf 1$ denotes the identify operator
and $b\geq m{\bf 1}$ then for any $1\leq s <2$
$$(b^{1/2}Tb^{1/2})^s\geq m^{s-1}b^{1/2}T^sb^{1/2}.$$
\end{lemma}

\begin{proof}
One can prove a weaker version of part (i) using singular values as a special
case of \cite{FK} Lemma 4.5. However, we feel that the stronger version has
some independent interest.
Now (i) is equivalent to:
$$\left(\left(\frac{b}{M}\right)^{1/2}T\left(\frac{b}{M}\right)^{1/2}\right)^s
\leq \left(\frac{b}{M}\right)^{1/2}T^s\left(\frac{b}{M}\right)^{1/2}.$$
So we can assume that $M=1$ and therefore $b \leq{\bf 1} $. Letting $A=b^{1/2}$ we
have $0\leq A\leq{\bf 1} $ and we want:
$$(ATA)^s\leq AT^sA.$$
Equivalently we want:
$$(ATA)(ATA)^{s-1}\leq ATT^{s-1}A$$
or, letting $r=s-1$ we want:
$$(ATA)(ATA)^r\leq ATT^rA$$
for $0\leq r <1.$
Using the integral formula for the $r^{th}$ power of a positive operator,
we want:
$$\int_0^\infty t^{-r}(ATA)(1+tATA)^{-1}ATAdt\leq
\int_0^\infty t^{-r}AT(1+tT)^{-1}TAdt$$
which would follow from:
$$\int_0^\infty t^{-r}[AT(A(1+tATA)^{-1}A-(1+tT)^{-1})TA]dt \leq 0.$$
So, it would be enough to see that:
$$A(1+tATA)^{-1}A\leq (1+tT)^{-1}.$$
Since the left hand side of this inequality is a norm-continuous function of
$A$, we can approximate $A$ by a sequence $\{A_n\}$ with $0<A_n\leq1$.
Then it suffices to prove that:
$$A_n(1+tA_nTA_n)^{-1}A_n\leq (1+tT)^{-1}.$$
Or,
$$(1+tA_nTA_n)\geq A_n(1+tT)A_n,$$
or,
$${\bf 1}\geq A_n^2.$$
So, (i) holds.

The argument for (ii) is very similar but easier. As in the proof of (i) we
can assume $m=1$ and letting $A=b^{1/2}$ we have $A\geq {\bf 1}$ and we want:
$$(ATA)^s\geq AT^sA$$
for $1\leq s<2$. We argue as above with all of the inequalities reversed.
Since $A\geq {\bf 1}$ it is invertible and we need no approximations. Our final
line for the argument then becomes ${\bf 1}\leq A^2$ and so (ii) is done.
\end{proof}

\begin{lemma}
For $T\geq 0$ in ${\mathcal L}^{(1,\infty)}$ and any
$b$ in $\mathcal N$ with $b\geq m{\bf 1}>0$,
$$\lim_{s\to 1^+}[(s-1)\tau(bT^s)-(s-1)\tau((b^{1/2}Tb^{1/2})^s)]=0$$
\end{lemma}

\begin{proof}
Let $M=||b||$ then by Lemma 3.3
$$(M^{s-1}-1)\tau(b^{1/2}T^sb^{1/2})\geq
\tau[(b^{1/2}Tb^{1/2})^s-b^{1/2}T^sb^{1/2}]
\geq (m^{s-1}-1)\tau[b^{1/2}T^sb^{1/2}].$$
Hence
$$(M^{s-1}-1)(s-1)\tau(b^{1/2}T^sb^{1/2})\geq
(s-1)\tau[(b^{1/2}Tb^{1/2})^s]-(s-1)\tau[b^{1/2}T^sb^{1/2}]$$
$$
\geq (m^{s-1}-1)(s-1)\tau(b^{1/2}T^sb^{1/2}).$$
Now let $s\to 1^+$:
$$0\geq \limsup_{s\to 1^+}\left((s-1)\tau[(b^{1/2}Tb^{1/2})^s]-(s-1)\tau[
bT^s]\right)$$
$$\geq\liminf_{s\to 1^+}\left((s-1)\tau[(b^{1/2}Tb^{1/2})^s]-
(s-1)\tau[bT^s]\right)\geq 0.$$
\end{proof}

\begin{lemma}
If $b\geq 0, T\geq 0$, $T\in{\mathcal L}^{(1,\infty)}$ and $b\in\mathcal N$
then there is a constant $C>0$ depending on $b,T$ such that for all
$0<\epsilon<1$.
$$ \limsup_{s\to 1^+}\left|(s-1)\tau[(b^{1/2}Tb^{1/2})^s]-(s-1)\tau
[((b+\epsilon)^{1/2}T(b+\epsilon)^{1/2})^s]\right|
\leq C\epsilon^{\frac{1}{4}}.$$
\end{lemma}

\begin{proof} To shorten the notation let
$A=b^{1/2}Tb^{1/2}$ and $B=(b+\epsilon)^{1/2}T(b+\epsilon)^{1/2})$
so that there is an $M>0$ such that $||A||_s\leq M||T||_s$
and  $||B||_s\leq M||T||_s$ for all $0<\epsilon<1$ and  $1<s<2$,
where $||.||_s$ is the Schatten class norm.
Then
$$\left|\tau[(b^{1/2}Tb^{1/2})^s]-\tau
[((b+\epsilon)^{1/2}T(b+\epsilon)^{1/2})^s]\right|\leq ||A^s-B^s||_1
$$
and
$$
||A^s-B^s||_1\leq ||A^{s/2}(A^{s/2}-B^{s/2})||_1 +||(A^{s/2}-B^{s/2})B^{s/2}||_1$$
Apply the \cite{BKS} inequality to the RHS of the previous line
(for a discussion of this inequality
for operator ideals in semifinite von Neumann algebras see
the references in \cite {CPS}) using $1>s/2$ to obtain
\begin{eqnarray*}
||A^s-B^s||_1&\leq&
||A^{s/2}||_2||A^{s/2}-B^{s/2}||_2 +||A^{s/2}-B^{s/2}||_2||B^{s/2}||_2\\
&\leq& ||A^{s/2}||_2|||A-B|^{s/2}||_2 +|||A-B|^{s/2}||_2||B^{s/2}||_2\\
&=&||A||_s^{s/2}||A-B||_s^{s/2}+||A-B||_s^{s/2}||B||_s^{s/2}\\
&\leq& 2M^{s/2}||T||_s^{s/2}||A-B||_s^{s/2}\\
&=& 2M^{s/2}(\tau(T^s))^{1/2}||A-B||_s^{s/2}.
\end{eqnarray*}
Hence
$$\left|(s-1)\tau(b^{1/2}Tb^{1/2})^s-(s-1)\tau
[((b+\epsilon)^{1/2}T^s(b+\epsilon)^{1/2})^s]\right|$$
$$\leq 2M^{s/2}((s-1)\tau(T^s))^{1/2}[(s-1)||A-B||_s^{s}]^{1/2}.$$
Now using the argument at the beginning of this section
there is a $K>0$ depending only on $b,T$
such that
$$ \limsup_{s\to 1^+} 2M^{s/2}((s-1)\tau(T^s))^{1/2}\leq K.$$
On the other hand
\begin{eqnarray*}
||A-B||_s&\leq& ||b^{1/2}T(b^{1/2}-(b+\epsilon)^{1/2})||_s
+||((b+\epsilon)^{1/2}-b^{1/2})T(b+\epsilon)^{1/2})||_s\\
&\leq& ||b^{1/2}||\ ||T||_s||b^{1/2}-(b+\epsilon)^{1/2}||
+||(b+\epsilon)^{1/2}-b^{1/2}||\ ||T||_s||(b+\epsilon)^{1/2}||\\
&\leq& K_2\sqrt\epsilon||T||_s=K_2\sqrt\epsilon(\tau(T^s))^\frac{1}{s}
\end{eqnarray*}
for some constant $K_2>0$.
Thus
$$
\limsup_{s\to 1^+}[(s-1)||A-B||_s^s]^\frac{1}{2}
\leq
\limsup_{s\to 1^+}[(s-1)\tau(T^s)]^\frac{1}{2}(K_2\sqrt\epsilon)^\frac{s}{2}
\leq (const)\epsilon^{1/4}.$$
as required.
\end{proof}

\begin{prop}
If $b\geq 0, T\geq 0$, $T\in{\mathcal L}^{(1,\infty)}$ and $b\in\mathcal N$
then
$\lim_{s\to 1^+}(s-1)\tau(bT^s)$
exists if and only if $\lim_{s\to 1^+}(s-1)\tau((b^{1/2}Tb^{1/2})^s)$
exists and in this case they are equal. Moreover, in any case
for any $\tilde\omega\in L^\infty(\IR)^*$ satisfying
the conditions (1), (2), (3), (4) of Corollary 1.6.
$$\tilde\omega-\lim_{r\to \infty}\frac{1}{r}\tau(bT^{1+\frac{1}{r}})=
\tilde\omega-\lim_{r\to\infty}
\frac{1}{r}\tau((b^{1/2}Tb^{1/2})^{1+\frac{1}{r}}).$$
\end{prop}

\begin{proof}
It suffices to prove:
$$\limsup_{r\to \infty}
\left|\frac{1}{r}\tau(bT^{1+\frac{1}{r}})
-\frac{1}{r}\tau((b^{1/2}Tb^{1/2})^{1+\frac{1}{r}})\right|=0.$$
Now,
\begin{eqnarray*}
&&\limsup_{r\to \infty}
\frac{1}{r}\left|\tau(bT^{1+\frac{1}{r}})
-\tau((b^{1/2}Tb^{1/2})^{1+\frac{1}{r}})\right|\\
&\leq&\limsup_{r\to \infty}\frac{1}{r}\left|\tau(bT^{1+\frac{1}{r}})
-\tau((b+\epsilon)T^{1+\frac{1}{r}})\right|\\
&\;&+\limsup_{r\to \infty}
\frac{1}{r}\left|
\tau((b+\epsilon)^{1/2}T^{1+\frac{1}{r}}(b+\epsilon)^{1/2})
-\tau[((b+\epsilon)^{1/2}T(b+\epsilon)^{1/2})^{1+\frac{1}{r}}]\right|\\
&\;&\;\;+\limsup_{r\to \infty}
\frac{1}{r}\left|
\tau[((b+\epsilon)^{1/2}T(b+\epsilon)^{1/2})^{1+\frac{1}{r}}]
-\tau((b^{1/2}Tb^{1/2})^{1+\frac{1}{r}})\right|\\
&\leq&\limsup_{r\to \infty}\frac{1}{r}\tau(T^{1+\frac{1}{r}})
\epsilon +0+C\epsilon^{1/4}
\end{eqnarray*}
by Lemmas 3.4 and 3.5. As this holds for all $\epsilon>0$
we are done.
\end{proof}

\begin{cor}
If $b\geq 0, T\geq 0$, $T\in{\mathcal L}^{(1,\infty)}$ and $b\in\mathcal N$
then if any one of the following limits exist
they all do and if $\omega$ is chosen to satisfy the conditions
of Theorem 1.5
they  are all equal to $\tau_\omega(bT)$\\
(1) $\lim_{t\to\infty}\frac{1}{\log(1+t)}\int_0^t\mu_s(b^{1/2}Tb^{1/2}) ds$\\
(2) $\lim_{r\to \infty}\frac{1}{r}\tau(bT^{1+\frac{1}{r}})$\\
(3) $\lim_{r\to \infty}\frac{1}{r}\tau((b^{1/2}Tb^{1/2})^{1+\frac{1}{r}})$
\end{cor}

\begin{proof}
The simultaneous existence and equality of (2) and (3)
 follows from Proposition 3.6. If
(3) exists then (1) exists and is equal to (3) by 
the second part of Theorem 3.1.

Conversely, if (1) exists then it equals $\tau_\omega(b^{1/2}Tb^{1/2})$
by definition. Then applying Lemma 3.2(i),
we have (1) equal to  $\tau_\omega(bT)$ and so for all $\epsilon>0$
there is an $M>0$ such that for
$t\geq M$
$$ \tau_\omega(bT)-\epsilon\leq
\frac{1}{\log(1+t)}\int_0^t\mu_s(b^{1/2}Tb^{1/2}) ds
\leq\tau_\omega(bT)+\epsilon.$$
Hence for $t\geq M$
$$ (\tau_\omega(bT)-\epsilon)\int_0^t\frac{1}{1+s}ds\leq
\int_0^t\mu_s(b^{1/2}Tb^{1/2}) ds\leq
(\tau_\omega(bT)+\epsilon)\int_0^t\frac{1}{1+s}ds.
$$
Following \cite{P} introduce three functions
$$
g_2(t)=\left\{\begin{array}{ll}
\frac{1}{M}\int_0^M\mu_s(b^{1/2}Tb^{1/2})ds & \;\;\;\;\mbox{if } t<M \\
\mu_t(b^{1/2}Tb^{1/2})  & \;\;\;\;\mbox{if } t\geq M
 \end{array}\right.$$
$$
g_1(t)=\left\{\begin{array}{ll}
(\tau_\omega(bT)-\epsilon)
\frac{1}{M}\int_0^M\frac{1}{1+s}ds & \mbox{if } t<M \\
(\tau_\omega(bT)-\epsilon)\frac{1}{1+t}  & \mbox{if } t\geq M
 \end{array}\right.$$
$$
g_3(t)=\left\{\begin{array}{ll}
(\tau_\omega(bT)+\epsilon)
\frac{1}{M}\int_0^M\frac{1}{1+s}ds & \mbox{if } t<M \\
(\tau_\omega(bT)+\epsilon)\frac{1}{1+t}  & \mbox{if } t\geq M
 \end{array}\right.$$
Then $g_1\prec\prec g_2\prec\prec g_3$
and thus
$g_1^{1+\frac{1}{r}}\prec\prec
g_2^{1+\frac{1}{r}}\prec\prec g_3^{1+\frac{1}{r}}$.
So we have for $t\geq M$
\begin{eqnarray*}
&&(\tau_\omega(bT)-\epsilon)^{1+\frac{1}{r}}
\left[M\left(\frac{1}{M}\int_0^M\frac{1}{1+s}ds\right)^{1+\frac{1}{r}}
+\int_M^t(\frac{1}{1+s})^{1+\frac{1}{r}}ds\right]\\
&\leq&M\left(\frac{1}{M}\int_0^M\mu_s(b^{1/2}Tb^{1/2})ds\right)^{1+\frac{1}{r}}
+\int_M^t \mu_s(b^{1/2}Tb^{1/2})^{1+\frac{1}{r}}ds\\
&\leq&(\tau_\omega(bT)+\epsilon)^{1+\frac{1}{r}}
\left[M\left(\frac{1}{M}\int_0^M\frac{1}{1+s}ds\right)^{1+\frac{1}{r}}
+\int_M^t(\frac{1}{1+s})^{1+\frac{1}{r}}ds\right].
\end{eqnarray*}
Let $t\to\infty$ so that
\begin{eqnarray*}
&&(\tau_\omega(bT)-\epsilon)^{1+\frac{1}{r}}
\left[M\left(\frac{1}{M}\int_0^M\frac{1}{1+s}ds\right)^{1+\frac{1}{r}}
+r(\frac{1}{1+M})^{\frac{1}{r}}\right]\\
&\leq&M\left(\frac{1}{M}\int_0^M\mu_s(b^{1/2}Tb^{1/2})ds\right)^{1+\frac{1}{r}}
+\tau((b^{1/2}Tb^{1/2})^{1+\frac{1}{r}})-
\int_0^M \mu_s(b^{1/2}Tb^{1/2})^{1+\frac{1}{r}}ds\\
&\leq&(\tau_\omega(bT)+\epsilon)^{1+\frac{1}{r}}
\left[M\left(\frac{1}{M}\int_0^M\frac{1}{1+s}ds\right)^{1+\frac{1}{r}}
+r(\frac{1}{1+M})^{\frac{1}{r}}\right].
\end{eqnarray*}
Multiply by $\frac{1}{r}$ and let $r\to \infty$,
$$\tau_\omega(bT)-\epsilon\leq
\lim_{r\to\infty}\frac{1}{r}\tau((b^{1/2}Tb^{1/2})^{1+\frac{1}{r}})
\leq(\tau_\omega(bT)+\epsilon).$$
Hence the result.
\end{proof}

\begin{thm}
Let $A\in\mathcal N$, $T\geq 0$, $T\in{\mathcal L}^{(1,\infty)}.$\\
(i) If $\lim_{s\to 1^+}(s-1)\tau(AT^s)$
exists then it is equal to $\tau_\omega(AT)$ where we choose $\omega$ as in
the proof of Theorem 3.1.\\
(ii) More generally, if we choose functionals $\omega$ and $\tilde\omega$ as in
the proof of Theorem 3.1 then
$$\tilde\omega-\lim_{r\to\infty}
\frac{1}{r}\tau(AT^{1+\frac{1}{r}})=\tau_\omega(AT).$$
\end{thm}

\begin{proof} For part (i) we first assume that $A$ is self adjoint.
Write $A=a^+-a^-$ where $a^\pm$ are positive.
Choose $\tilde\omega$ as in the proof of
Theorem 3.1, then
\begin{eqnarray*}
\lim_{s\to 1^+}(s-1)\tau(AT^s)&=&
\tilde\omega-\lim_{r\to\infty}\frac{1}{r}
\tau(AT^{1+\frac{1}{r}})\\
&=&\tilde\omega-\lim_{r\to\infty}\frac{1}{r}\tau(a^+T^{1+\frac{1}{r}})
-\tilde\omega-\lim_{r\to\infty}\frac{1}{r}\tau(a^-T^{1+\frac{1}{r}})\\
&=& \tau_\omega(a^+T)-\tau_\omega(a^-T)\\
&=&\tau_\omega(AT).
\end{eqnarray*}
Here the third equality uses first Proposition 3.6 and then Theorem 3.1.
The reduction from the general case to the self-adjoint case now follows
in a similar way.

For part (ii), we assume that $A$ is positive. By Lemma 3.2(i), Theorem 3.1,
and Proposition 3.6 we have
\begin{eqnarray*}
\tau_\omega(AT)&=&\tau_\omega(A^{1/2}TA^{1/2})=\tilde\omega-
\lim_{r\to \infty}\frac{1}{r}\tau((A^{1/2}TA^{1/2})^{1+\frac{1}{r}})\\
&=&\tilde\omega-\lim_{r\to \infty}\frac{1}{r}\tau(AT^{1+\frac{1}{r}}).
\end{eqnarray*}
For general $A$ we reduce to the case $A$ positive as in the proof of part (i).
\end{proof}

\section{\bf The heat semigroup formula}

Throughout this section  $T\geq 0$.
We define $e^{-T^{-2}}$ as the operator that is zero on $\ker T$
and on $\ker T^\perp$ is defined in the usual way by the functional
calculus. We remark that if $T\geq 0$, $T\in{\mathcal L}^{(p,\infty)}$
for some $p\geq 1$ then  $e^{-tT^{-2}}$ is trace class for all $t>0$.

Our aim in this section is to prove the following

\begin{thm}
If $A\in\mathcal N$, $T\geq 0$, $T\in{\mathcal L}^{(1,\infty)}$
then,
$$\omega-\lim_{\lambda\to\infty}\lambda^{-1}\tau(Ae^{-\lambda^{-2}T^{-2}})=
\Gamma(3/2)\tau_\omega(AT)$$
for $\omega\in L^\infty(\IR^*_+)^*$ satisfying the conditions of Theorem 1.5.
\end{thm}

Let $\zeta_A(p+\frac{1}{r})=\tau(AT^{p+\frac{1}{r}})$.
Notice that $\frac{1}{2}\Gamma(\frac{p}{2})
\tilde\omega-\lim_{r\to\infty} \frac{1}{r}\zeta_A(p+\frac{1}{r})$
always exists. Hence we can reduce the hard part of the
proof of Theorem 4.1 to the following preliminary result.

\begin{prop}
If $A\in\mathcal N$, $T\geq 0$, $T\in{\mathcal L}^{(p,\infty)}$, $1\leq p<\infty$
then, choosing $\omega$ and $\tilde\omega$
as in the proof of theorem 3.1, we have
$$\omega-\lim \frac{1}{\lambda}\tau(Ae^{-T^{-2}\lambda^{-2/p}})
=\frac{1}{2}\Gamma(\frac{p}{2})
\tilde\omega-\lim_{r\to\infty} \frac{1}{r}\zeta_A(p+\frac{1}{r}).$$
\end{prop}

\begin{proof}
We have, using the Laplace transform,
$$T^s = \frac{1}{\Gamma(s/2)}\int_0^\infty t^{s/2 -1}e^{-tT^{-2}}dt.$$
Then
$$\zeta_A(s)=\tau(AT^s)=
\frac{1}{\Gamma(s/2)}\int_0^\infty t^{s/2 -1}\tau(Ae^{-tT^{-2}})dt.$$
Make the change of variable $t=1/\lambda^{2/p}$ so that the preceding
formula becomes
$$\frac{p}{2}\Gamma(s/2)\zeta_A(s)=
\int_0^\infty \lambda^{-\frac{s}{p}-1}
\tau(Ae^{-\lambda^{-2/p}T^{-2}})d\lambda.$$
We split this integral into two parts,
$\int_0^1$ and $\int_1^\infty$
and call the first integral $R(r)$ where $s=p+\frac{1}{r}$.
Then
$$R(r)=\int_0^1 \lambda^{-\frac{1}{pr}-2}
\tau(Ae^{-\lambda^{-2/p}T^{-2}})d\lambda=\int_1^\infty
t^{\frac{p}{2}+\frac{1}{2r} -1}\tau(Ae^{-tT^{-2}})dt.$$
The integrand decays exponentially in $t$ as $t\to\infty$
because 
$T^{-2}\geq {\Vert T^2\Vert}^{-1}\bf 1$
so that
$$\tau(Ae^{-tT^{-2}})\leq \tau(Ae^{-T^{-2}}e^{-\frac{t-1}{\Vert T^{2}\Vert}}).
$$
Then we can conclude that $R(r)$ is bounded independently of $r$ and so
$\lim_{r\to\infty}\frac{1}{r}R(r)=0$.
For the other integral
the change of variable
$\lambda = e^\mu$ gives
$$\int_1^\infty \lambda^{-\frac{1}{pr}-2}
\tau(Ae^{-\lambda^{-2/p}T^{-2}})d\lambda
=\int_0^\infty e^{-\frac{\mu}{pr}} d\beta(\mu)$$
where
$\beta(\mu)=\int_0^\mu e^{-v}\tau(Ae^{-e^{-\frac{2}{p}v}T^{-2}})dv$.
Hence we can now write
$$\frac{p}{2}\Gamma((p+\frac{1}{r})/2)\zeta_A(p+\frac{1}{r})
=\int_0^\infty e^{-\frac{\mu}{pr}} d\beta(\mu)+R(r).$$
Now consider
$$\frac{p}{2}\tilde\omega-\lim_{r\to\infty}\frac{1}{r}
\Gamma(\frac{p}{2}+\frac{p}{2r})\zeta_A(p+\frac{1}{r})=
\frac{p}{2}
\Gamma(p/2)\tilde\omega-\lim_{r\to\infty}\frac{1}{r}\zeta_A(p+\frac{1}{r}).$$
Then
$$\frac{p}{2}
\Gamma(p/2)\tilde\omega-\lim_{r\to\infty}\frac{1}{r}\zeta_A(p+\frac{1}{r})
=p\tilde\omega-\lim_{r\to\infty}\frac{1}{pr}
\int_0^\infty e^{-\mu/pr}d\beta(\mu)$$
(remembering that the term $\frac{1}{r}R(r)$ has limit
zero as $r\to \infty$).
By dilation invariance and
Theorem 2.2 we then have
$$\frac{p}{2}
\Gamma(p/2)\tilde\omega-\lim_{r\to\infty}\frac{1}{r}
\zeta_A(p+\frac{1}{r})=p\tilde\omega-\lim_{\mu\to\infty}
 \frac{\beta(\mu)}{\mu}. \eqno (4.0)$$
Making the change of variable $\lambda = e^v$
in the expression for $\beta(\mu)$
we get
$$\frac{\beta(\mu)}{\mu}=\frac{1}{\mu}\int_1^{e{^\mu}}
\lambda^{-2}\tau(Ae^{-T^{-2}\lambda^{-2/p}})d\lambda$$
Make the substitution $\mu =\log t$ so the RHS becomes
$$\frac{1}{\log
t}\int_1^t\lambda^{-2}\tau(Ae^{-T^{-2}\lambda^{-2/p}})d\lambda=g_1(t)$$
This is the Cesaro mean of
$$g_2(\lambda)=\frac{1}{\lambda}\tau(Ae^{-T^{-2}\lambda^{-2/p}}).$$
So as we chose $\omega\in L^\infty(\IR^*_+)^*$ to be $G_2$ and $M$ invariant
we have $\omega(g_1)=\omega(g_2)$. Recalling that we choose
$\tilde\omega$ to be related to $\omega$ as in Theorem 3.1
and so using (4.0) we obtain
$$\omega-\lim \frac{1}{\lambda}\tau(Ae^{-T^{-2}\lambda^{-2/p}})
=\frac{1}{2}\Gamma(\frac{p}{2})
\tilde\omega-\lim_{r\to\infty} \frac{1}{r}\zeta_A(p+\frac{1}{r}).$$
\end{proof}

To prove the theorem consider first the
case where $A$ is bounded, $A\geq 0$ and use
the Proposition 4.2 and Theorem 3.8 to assert that
$$\Gamma(3/2)\tau_\omega(AT)=
\Gamma(3/2)\tilde\omega-\lim_{r\to\infty}\frac{1}{r}\tau(AT^{1+\frac{1}{r}})
=\omega-\lim_{\lambda\to\infty}\lambda^{-1}
\tau(Ae^{-\lambda^{-2}T^{-2}}).
$$
Then for self adjoint $A$
write $A=a^+-a^-$ where $a^\pm$ are positive so that
$$\Gamma(3/2)\tau_\omega(AT)=
\Gamma(3/2)(\tau_\omega(a^+T)-\tau_\omega(a^-T))$$
$$=
\omega-\lim_{\lambda\to\infty}\lambda^{-1}
\tau(a^+e^{-\lambda^{-2}T^{-2}})-
\omega-\lim_{\lambda\to\infty}\lambda^{-1}\tau(a^-e^{-\lambda^{-2}T^{-2}})
$$
$$=
\omega-\lim_{\lambda\to\infty}\lambda^{-1}\tau(Ae^{-\lambda^{-2}T^{-2}}).$$
We can extend to general bounded $A$ by a similar argument.

\subsection{\bf The `smaller' ideal}

The curious feature of our proof of this heat kernel formula of Connes
for the Dixmier trace is that we need to go via the zeta function
and hence need the pair of functionals $\tilde\omega$ and $\omega$
as in Theorem 3.1. There is a special case of the previous result
for which we can avoid the introduction of these functionals
and hence avoid using the full strength of the assumptions in
Theorem 1.5.

The operators $T\in{\mathcal L}^{(1,\infty)}$
satisfying $\mu_s(T)\leq C/s$ for some $C>0$
form an ideal as well. For this `smaller ideal',
which is the one that usually arises in geometric
applications,
there is a direct proof of
a special case of the heat kernel
formula which does not use the zeta function.

For simplicity we restrict to $A=1$.
This direct proof uses the
 Laplace transform:
$T=\frac{1}{\Gamma(1/2)}\int_0^\infty u^{-1/2}
e^{-uT^{-2}}du$ (with our usual convention
that $e^{-T^{-2}}$ is defined to be zero on $\ker T$).
Thus we have
$$\frac{\Gamma(3/2)}{\log(1+t)}\int_0^t \mu_s(T)ds
=\frac{1}{2\log(1+t)}\int_0^t\int_0^\infty u^{-1/2}
e^{-u/\mu_s(T^2)}duds \eqno(4.1)$$

Using the basic fact that if $f$ is increasing
$\mu_s(f(T))=f(\mu_s(T))$ \cite{FK}
we have
$$\frac{1}{\log(1+t)}\int_0^t\lambda^{-2}\tau(e^{-\lambda^{-2}T^{-2}})
d\lambda=\frac{1}{\log(1+t)}\int_0^t\int_0^\infty
\lambda^{-2}e^{-\lambda^{-2}/\mu_s(T^2)}dsd\lambda$$
and we have to show that this has the same $\omega$ limit as (4.1).
Change variable in this integral by $u=\lambda^{-2}$ then
$$\frac{1}{\log(1+t)}\int_0^t\lambda^{-2}\tau(e^{-\lambda^{-2}T^{-2}})
d\lambda=\frac{1}{2\log(1+t)}\int_{1/t^2}^\infty\int_0^\infty u^{-1/2}
e^{-u/\mu_s(T^2)}dsdu. \eqno(4.2)$$

Subtract (4.1) and (4.2)
and rewrite the difference as
$$\frac{1}{2\log(1+t)}\int_0^\infty\int_0^\infty
(-\chi_{[0,t]}(s)\chi_{[0,1/t^2]}(u)
+\chi_{[1/t^2,\infty)}(u)\chi_{[t,\infty)}(s)) u^{-1/2}
e^{-u/\mu_s(T^2)}duds. \eqno(4.3)$$

To prove equality of the $\omega$-limits of (4.1) and (4.2)
we have to estimate the two integrals in (4.3).
The first of these is
$$\frac{1}{2\log(1+t)}\int_0^{1/t^2}\int_0^{t} u^{-1/2}
e^{-u/\mu_s(T^2)}duds.$$
As $\mu_s(T^2)\to 0$ as $s\to \infty$
we can assume there is a constant $C$
such that $e^{-u/\mu_s(T^2)}\leq Ce^{-u}$.
Thus the integral is bounded by
$$\frac{1}{2\log(1+t)}\int_0^{t}(\int_0^{1/t^2} u^{-1/2} Ce^{-u}du)ds
=\frac{1}{2\log(1+t)}t\int_0^{1/t^2} u^{-1/2} Ce^{-u}du.$$
Now
$$\gamma(\frac{1}{2},\frac{1}{t^2})=\int_0^{1/t^2} u^{-1/2} e^{-u}du$$
is the incomplete $\Gamma$ function which has an expansion
of the form (see \cite{AS})
$$\gamma(\frac{1}{2},\frac{1}{t^2})
=\frac{1}{t}\sum_0^\infty\frac{(-1)^n}{n!}\frac{1}{t^{2n}(\frac{1}{2}+n)}.$$
So we conclude that
$$t\int_0^{1/t^2} u^{-1/2} Ce^{-u}duds
=C\sum_0^\infty\frac{(-1)^n}{n!}\frac{1}{t^{2n}(\frac{1}{2}+n)}$$
which is bounded as $t\to\infty$. Thus as $t\to \infty$
$$\frac{1}{2\log(1+t)}\int_0^{1/t^2}\int_0^{t} u^{-1/2}
e^{-u/\mu_s(T^2)}duds\to 0.$$

For the second integral
in (4.3) we first make a number of preliminary
observations. We make some changes of variable in
letting $r=s/t$ and $v=ut^2$.
Then we find that
$$\int_{1/t^2}^\infty\int_{t}^\infty u^{-1/2}
e^{-u/\mu_s(T^2)}dsdu
=\int_1^\infty\int_1^\infty v^{-1/2} e^{-v/t^2\mu_{rt}(T^2)} drdv.$$
Now we exploit the assumption that
$\mu_s(T)=O(1/s)$ and use $v^{-1/2}<1$.
Thus
$\mu_{rt}(T^2)\leq C/(rt)^2$ for some constant $C$ and
$$\int_1^\infty\int_1^\infty v^{-1/2} e^{-v/t^2\mu_{rt}(T^2)} drdv
\leq \int_1^\infty\int_1^\infty e^{-v r^2/C}dv dr$$
$$=\int_1^\infty \frac{1}{r^2}e^{-r^2/C} dr<\infty.$$

Dividing by $\log (1+t)$ and taking $t\to \infty$
shows that the second integral in (4.3)
gives a function of $t$
which vanishes at infinity.

Now choose $\omega\in L^\infty(\IR^*_+)^*$ satisfying
conditions (1), (2), (3), (6) of Theorem 1.5.
Taking the $\omega$-limit on (4.3) gives zero.
Writing $\tau_\omega(T)=\omega-\lim_{t\to\infty}\frac{1}{\log(1+t)}
\int_0^t\mu_s(T) ds$ we obtain, using the same reasoning as at the
end of Proposition 4.2,
the result that
$$\omega-\lim_{\lambda\to\infty}\lambda^{-1}\tau(e^{-\lambda^{-2}T^{-2}})=
\Gamma(3/2)\tau_\omega(T).$$

\section{\bf The ${\mathcal L}^{(p,\infty)}$-summable case}

If $T\in{\mathcal L}^{(p,\infty)}$ for $p>1$, $T\geq 0$
then $\mu_s(T)=O(\frac{1}{s^{1/p}})$.
Moreover $\tau(T^{p+\frac{1}{r}}) =\int_0^1 \lambda^{p+1/r} dN_T(\lambda)$
where $N_T(\lambda)=\tau(E(\lambda)-1)$ for $\lambda>0$ where 
$T=\int \lambda dE(\lambda)$ is the spectral resolution for $T$.

We now establish some  ${\mathcal L}^{(p,\infty)}$
versions of our previous results.

\begin{lemma} For $T\in{\mathcal L}^{(p,\infty)}$ and
$\omega$ and $\tilde\omega$ as in the proof of theorem 3.1
we have
$$p\tau_\omega(T^p)=\tilde\omega-\lim_{r\to\infty}
\frac{1}{r} \tau(T^{p+\frac{1}{r}}).$$
\end{lemma}

\begin{proof}
Set $\lambda=e^{-u/p}$ so that
$$\frac{1}{r} \tau(T^{p+\frac{1}{r}}) =p\frac{1}{pr}\int_0^\infty
e^{-u/rp}d\beta(u)$$
where $\beta(u)=\int_0^ue^{-v}dN_T(e^{-v/p})$.
So using dilation invariance:
$$\tilde\omega-\lim_{r\to\infty}\frac{1}{r} \tau(T^{p+\frac{1}{r}})
=p\tilde\omega-\lim_{r\to\infty}
\frac{1}{pr}\int_0^1
e^{-u/pr}d\beta(u)=p\tilde\omega-\lim_{u\to\infty}
\frac{\beta(u)}{u}$$
by the weak$^*$-Karamata theorem.
Reasoning as in the proof of Theorem 3.1  and substituting
$\lambda=e^{-v/p}$ and
$u=\log t$ we have
$$\tilde\omega-\lim_{u\to\infty}
\frac{\beta(u)}{u}=\omega-\lim_{t\to\infty}\frac{1}{\log{t}}
\int_{t^{-1/p}}^1 \lambda^p dN_T(\lambda)$$
$$=\omega-\lim_{t\to\infty}\frac{1}{\log{t}}\tau(\chi_{(\frac{1}{t},\infty)}
(T^p)T^p)=\tau_\omega(T^p).$$
\end{proof}

\begin{cor}
Let $T\geq 0$, $T\in{\mathcal L}^{(p,\infty)}$ then
$$\omega-\lim \frac{1}{\lambda}\tau(e^{-T^{-2}\lambda^{-2/p}})
=\Gamma(1+\frac{p}{2})
\tilde\omega-\lim_{r\to\infty} \frac{1}{r}\tau(T^{p+\frac{1}{r}})$$
with the usual convention that $e^{-T^{-2}}$ is zero on $\ker T$.
\end{cor}

\begin{proof}
Combine Proposition 4.2 and Lemma 5.1.
\end{proof}

Our aim is now to prove
the ${\mathcal L}^{(p,\infty)}$ version of
Theorem 3.8 and the following result of Connes'.

\begin{thm}
If $A$ is bounded, $T\geq 0$, $T\in{\mathcal L}^{(p,\infty)}$
for $p\geq 1$
$$\omega-\lim_{\lambda\to\infty}\lambda^{-1}\tau(Ae^{-\lambda^{-2/p}T^{-2}})=
\Gamma(1+p/2)\tau_\omega(AT^p)$$
where $e^{-T^{-2}}$ is defined to be zero on $\ker T$.
\end{thm}

To this end let us consider the steps in the proof of Theorem
3.8. The key results are
Proposition 3.6 and Corollary 3.7.
Proposition 3.6 rests on the preceding lemmas.
These lemmas have analogues in the case of ${\mathcal L}^{(p,\infty)}$.
The first non-obvious extension is Lemma 3.3 which we replace by

\begin{lemma}
Let $0\leq T\in{\mathcal L}^{(1,\infty)}$
 and let $0\leq b\in {\mathcal N}$, $||b||\leq M$.

\noindent(i) For any $s\geq 1$
$$\mu_t(b^{1/2}Tb^{1/2})^s\leq M^{s-1} \mu_t(b^{1/2}T^sb^{1/2}),\quad t>0.$$
\noindent(ii)
If $b\ge m1$, then
$$\mu_t(b^{1/2}Tb^{1/2})^s\geq m^{s-1} \mu_t(b^{1/2}T^sb^{1/2}),\quad t>0.$$

\end{lemma}

\begin{proof}
The first result is a special case of [FK] Lemma 4.5.
To obtain the second result, we shall (without loss of generality)
 assume that $\Vert T\Vert \leq 1$.
Let $T=\int _0^1\lambda dE^T(\lambda)$ be the spectral decomposition of $T$.
Note that it follows from the assumption $T\in{\mathcal L}^{(1,\infty)}$ that
$\tau (E^T(1/n,1])<\infty$ for all $n\in {\mathbb N}$. We set for brevity
$$
p_n:=E^T(1/n,1],\quad q_n:=l(b^{-1/2}p_n)\bigvee r(p_nb^{-1/2}),\quad
{\mathcal N}_n:=q_n{\mathcal N}q_n,\quad n\in {\mathbb N},
$$
where $l(\cdot)$ and $r(\cdot)$ are left and right support projections
respectively. Note that ${\mathcal N}_n$ is a finite von Neumann algebra
and that restriction  of the trace $\tau$ on ${\mathcal N}_n$  is
semifinite for every  $n\in {\mathbb N}$. From assertion $(i)$ we have
$$
\mu_t(b^{-1/2}(p_nTp_n)^{-1}b^{-1/2})^s\leq
m^{-(s-1)}\mu_t(b^{-1/2}(p_nTp_n)^{-s}b^{-1/2}),\quad n\in {\mathbb N}.
\eqno (5.1)
$$
Note that $b^{-1/2}\ge \frac1M$ and therefore $b^{-1/2}(p_nTp_n)^{-1}b^{-1/2},
\ b^{-1/2}(p_nTp_n)^{-s}b^{-1/2}$ are invertible elements for all $n\ge 1$.

 Now we need a following simple observation: if $({\mathcal M},\tau)$ is a
finite von Neumann algebra and $0\leq x$ is an invertible
$\tau$-measurable operator affiliated with ${\mathcal M}$, then the elements
$x^{-1}$ and $\mu_{(\cdot)}(x)^{-1}$ are equimeasurable, or equivalently,
$\mu_{(\cdot)}(x^{-1})$ is the decreasing rearrangement of the function
$\mu_{(\cdot)}(x)^{-1}$. To see the validity of this observation, set for
brevity $f(\lambda):=\frac 1\lambda$, $x=\int_0^\infty\lambda dE^x(\lambda),\
y=f(x)=\int_0^\infty f(\lambda)dE^x(\lambda) =\int_0^\infty\lambda dE^y(\lambda)
$
and note that $E^y(\Delta)=E^x(f^{-1}(\Delta))$ for every Borel subset
$\Delta\subseteq [0,\infty)$. In particular,
$$
E^y(s,\infty)=E^x(f^{-1}(s,\infty))=E^x(0,\frac 1s)=1-E^x[\frac 1s, \infty),
\quad s>0,
$$
whence
$$
\lambda _s(y)=\tau(1)-\lambda_{\frac 1s -0}(x), \ s>0.
$$
If instead of
the algebra  $({\mathcal M},\tau)$ and the element $x$ we consider the von
Neumann algebra $L_\infty (0, \tau(1))$ and the element
$\mu_{(\cdot)}(x)$, then
the preceding equality becomes
$$
\lambda _s((\mu_{(\cdot)}(x))^{-1})=\tau(1)-\lambda_{\frac 1s -0}
(\mu_{(\cdot)}(x)), \ s>0,
$$
(where we use the notation $\lambda_{(\cdot )}$ for the classical distribution
function of the elements $(\mu_{(\cdot)}(x))^{-1}$ and $\mu_{(\cdot)}(x)$).
Our observation now follows from comparison of the two preceding equalities,
taking into account a crucial fact, namely that
$\lambda_{\frac 1s -0}(x)=\lambda_{\frac 1s -0}(\mu_{(\cdot)}(x))$
for all $s>0$. This latter fact easily follows from
the equality $\lambda_{s}(x)= \lambda_{s}(\mu_{(\cdot)} (x))$ and
the assumption that ${\mathcal M}$ is finite.

Now we can continue the proof of the lemma. From the inequality $(5.1)$
taking the inverses we get
$$
\mu_t^{-1}(b^{-1/2}(p_nTp_n)^{-1}b^{-1/2})^s\ge
m^{(s-1)}\mu_t^{-1}(b^{-1/2}(p_nTp_n)^{-s}b^{-1/2}),\quad t>0,\quad n\in
{\mathbb N},\quad s\geq 1.
$$
Since $0\leq x\leq y$ implies $\mu_{(\cdot)} (x)\leq \mu_{(\cdot)} (y)$ we
immediately infer from the preceding inequality
$$
\mu_{(\cdot)} \left (\mu_t^{-1}(b^{-1/2}(p_nTp_n)^{-1}b^{-1/2})^s\right) \ge
m^{(s-1)}\mu_{(\cdot)} \left ( \mu_t^{-1}(b^{-1/2}(p_nTp_n)^{-s}b^{-1/2})\right),
 \quad n\in {\mathbb N}.
$$
The elements
$(b^{-1/2}(p_nTp_n)^{-1}b^{-1/2})^s$ and $b^{-1/2}(p_nTp_n)^{-s}b^{-1/2}$
are invertible positive elements from ${\mathcal N}_n$, and by the
preceding observation we know that the elements
$\mu_{(\cdot)}^{-1}(b^{-1/2}(p_nTp_n)^{-1}b^{-1/2})^s$ and
$(b^{1/2}(p_nTp_n)b^{1/2})^s$
 (respectively, $ \mu_{(.)}^{-1}(b^{-1/2}(p_nTp_n)^{-s}b^{-1/2})$
 and $ b^{1/2}(p_nTp_n)^sb^{1/2})$) are equimeasurable, thus
the preceding inequality may be equivalently rewritten as
$$
\mu_{(\cdot)} \left ((b^{1/2}(p_nTp_n)b^{1/2})^s\right) \ge
m^{(s-1)}\mu_{(\cdot)} \left (b^{1/2}(p_nTp_n)^{s}b^{1/2}\right),
 \quad n\in {\mathbb N}.
$$
To complete the proof of the lemma it is sufficient to show that
$$
\mu_{(\cdot)} \left ((b^{1/2}(p_nTp_n)b^{1/2})^s\right) \to
\mu_{(\cdot)} \left ((b^{1/2}Tb^{1/2})^s\right)\eqno (5.2)
$$
and
$$
\mu_{(\cdot)} \left (b^{1/2}(p_nTp_n)^{s}b^{1/2}\right)  \to
\mu_{(\cdot)} \left (b^{1/2}T^{s}b^{1/2}\right)\eqno (5.3)
$$
in measure. Since $\mu_{(\cdot)}(x^s)=\mu_{(\cdot)}^s(x)$ for all
$x\in {\mathcal N}$ and all $s>0$, to establish the first convergence,
it is sufficient to show that
$$
\mu_{(\cdot)} \left (b^{1/2}(p_nTp_n)b^{1/2}\right) \to
\mu_{(\cdot)} \left (b^{1/2}Tb^{1/2}\right).
$$
To this end we shall need the following result (\cite{CS} Corollary 2.3).

If $E({\mathcal N})$ is a symmetric operator space associated with
a separable symmetric operator space $E(0,\infty)$, then
$\Vert xe_n\Vert _{E({\mathcal N})}\to 0$ and
$\Vert e_nx\Vert _{E({\mathcal N})}\to 0$
for every $x\in E({\mathcal N})$ and every sequence $\{e_n\}$ of orthogonal
projections in ${\mathcal N}$ decreasing to $0$.

 Consider the symmetric function space $L_1+L_\infty (0,\infty)$ and let $E$ be
its closed separable symmetric subspace obtained by taking the norm closure of
$L_1\cap L_\infty (0,\infty)$. It is easy to see that $E$ is a separable
symmetric function space (in a sense it is an analogue of the space $c_0$ of
all bounded sequences converging to $0$). It is clear from the cited result
from \cite{CS}
 and the definition of $p_n$ that
$\Vert T-Tp_n\Vert_{E({\mathcal N})}\to 0$ and
$\Vert T-p_nT\Vert_{E({\mathcal N})}\to 0$, whence
 $\Vert T-p_nTp_n\Vert_{E({\mathcal N})}\to 0$ and also
$$
\Vert b^{1/2}(p_nTp_n)b^{1/2}-b^{1/2}Tb^{1/2}\Vert_{E({\mathcal N})}\to 0.
$$
Using the continuity of embedding of any $E({\mathcal N})$ into the space
$\widetilde{\mathcal N}$ of all $\tau$-measurable operators affiliated with
${\mathcal N}$ we get from the preceding convergence that
$$
 b^{1/2}(p_nTp_n)b^{1/2}-b^{1/2}Tb^{1/2}\to 0
$$
in measure. Now using [FK], Lemma 3.4 (ii) and the fact
$$
\lim_{t\to \infty}\mu_{t}(b^{1/2}(p_nTp_n)b^{1/2})=
\lim_{t\to \infty} \mu_{t}(b^{1/2}Tb^{1/2})=0
$$
we get
$$
 \mu_{(\cdot )}(b^{1/2}(p_nTp_n)b^{1/2})-\mu_{(\cdot )}(b^{1/2}Tb^{1/2})\to 0
$$
in measure, i.e. $(5.2)$ is established. The proof of $(5.3)$ is very similar,
after we note that $(p_nTp_n)^{s}=(p_nT^sp_n)$ and therefore we omit
the details.
\end{proof}

Next, some remarks are needed for Lemma 3.5.
For the ${\mathcal L}^{(p,\infty)}$ case the statement reads
if $b\geq 0, T\geq 0$, $T\in{\mathcal L}^{(1,\infty)}$ with $b$ bounded
then there is a constant $C>0$ depending on $b,T$ such that for all
$0<\epsilon<1$.
$$ \limsup_{s\to p^+}\{(s-p)\tau(b^{1/2}Tb^{1/2})^s-(s-p)\tau
[((b+\epsilon)^{1/2}T(b+\epsilon)^{1/2})^s]\}
\leq C\epsilon^\frac{1}{4}.$$
For the proof we use the same argument for all $1<p<2$
but for $p\geq 2$ we use Cauchy-Schwartz in place of the BKS inequality
so that in fact the proof is more elementary.
The proofs of Proposition 3.6 and Corollary 3.7
also generalise to give
us the

\begin{prop}
If $b\geq 0, T\geq 0$, $T\in{\mathcal L}^{(p,\infty)}$ with $b$ bounded
then
$\lim_{s\to p^+}(s-p)\tau(bT^s)$
exists if and only if $\lim_{s\to p^+}(s-p)\tau((b^{1/2}Tb^{1/2})^s)$
exists and in this case they are equal. Moreover, in any case
for any $\omega\in L^\infty(\IR_+)^*$ chosen to satisfy the conditions of
theorem 1.2
$$\omega-\lim_{r\to \infty}\frac{1}{r}\tau(bT^{p+\frac{1}{r}})=
\omega-\lim_{r\to \infty}\frac{1}{r}\tau((b^{1/2}Tb^{1/2})^{p+\frac{1}{r}}).$$
\end{prop}

Now the proof of Theorem 3.8(i)
generalises to give the

\begin{thm}
If $A$ is bounded, $T\geq 0$, $T\in{\mathcal L}^{(p,\infty)}$
and
$$\lim_{s\to p^+}(s-p)\tau(AT^s)$$
exists then it is equal to $p\tau_\omega(AT^p)$.
\end{thm}

Finally it is now straightforward to extend the arguments
we used in the proof of Theorem 4.1 to
prove Theorem 5.3.

\section{\bf Application to spectral flow and index formulae}

We fix an unbounded self-adjoint operator $D_0$
on $H$ affiliated with
$\mathcal N$.
Recalling the discussion in the introduction we have:

\begin{defn*}
 We say
that $(\mathcal N,D_0)$ is an {\bf odd unbounded
${\mathcal L}^{(1,\infty)}$-summable
Breuer-Fredholm module}
for a Banach $*$-algebra $\mathcal A$ if $\mathcal A$ is represented in
$\mathcal N$ and if $(1+D^2_0)^{-1/2}\in {\mathcal L}^{(1,\infty)}$
and $[D_0,a]$ is bounded for all $a$ in a dense $*$-subalgebra of $\mathcal A$.
\end{defn*}

Recall that these assumptions imply that $a$ leaves the domain of $D_0$
invariant.
In this section we apply our results to ${\mathcal L}^{(1,\infty)}$
summable Breuer-Fredholm modules in order to establish
a relationship between the formula for spectral flow
in \cite{CP2} and the
formula in \cite{CM}. In \cite{CM} assumptions
are made about the discreteness of the spectrum of $D_0$
which are clearly unrealistic when $\mathcal N$
is not type I.

We now summarise some well known notions ({\it cf} \cite{PR}).
Let ${\mathcal K}_{\mathcal N }$ be the 
$\tau$-compact operators in ${\mathcal N}$
(that is the norm closed ideal generated by the projections
$E\in\mathcal N$ with $\tau(E)<\infty$)
and
$\pi:{\mathcal N}\to {\mathcal N}/{\mathcal K_{\cn}}$
the canonical mapping. A Breuer-Fredholm operator is one that
maps to an invertible operator under $\pi$.
For a unitary $u\in \mathcal A$ the path
$$D_t^u:=(1-t)\,D_0+tuD_0u^*$$
of unbounded self-adjoint
Breuer-Fredholm operators is continuous in the sense that
$$F_t^u:=D_t^u\left(1+(D_t^u)^2\right)^{-\frac{1}{2}}$$
is a continuous path of self-adjoint Breuer-Fredholm operators in $
\mathcal N$.
Recall that the Breuer-Fredholm index of
a Breuer-Fredholm operator $F$ is defined by
$$ind(F)=\tau(Q_{ker F})-\tau(Q_{coker F})$$
where $Q_{ker F}$ and $Q_{coker F}$ are the projections onto
the kernel and cokernel of $F$.

\begin{defn*}
If $\{F_t\}$ is a continuous path of self-adjoint Breuer-Fredholm operators
in $\mathcal N$, then the definition of the {\bf spectral flow} of the
path, $sf(\{F_t\})$ is based on the following sequence of observations
in \cite{P1}:

\noindent 1. The map $t\mapsto sign(F_t)$ is usually discontinuous as is the
projection-valued mapping $t\mapsto P_t=\frac{1}{2}(sign(F_t)+1).$

\noindent 2. However,  $t\mapsto \pi(P_t)$ is continuous.

\noindent 3. If $P$ and $Q$ are projections in $\mathcal N$ and
$||\pi(P)-\pi(Q)||<1$
then $$PQ:rng(Q)\to rng(P)$$ is a Breuer-Fredholm operator and so
$ind(PQ)\in \real$ is well-defined.

\noindent 4. If we partition the parameter interval of $\{F_t\}$ so
that the $\pi(P_t)$ do not vary much in norm on each subinterval
of the partition then $$sf(\{F_t\}):=\sum_{i=1}^n ind(P_{t_{i-1}}P_{t_i})$$
is a well-defined and (path-) homotopy-invariant number which
agrees with the usual notion of spectral flow in the type $I_\infty$
case.
\end{defn*}

We denote by $sf(D_0,uD_0u^*)=sf(\{F_t\})$ the {\it spectral flow} of this path
[P1;P2] which is an integer in the
${\mathcal N}={\mathcal B}({\mathcal H})$ case and a real
number in the general semifinite case. This real number $sf(D_0,uD_0u^*)$
recovers the pairing of the $K$-homology class
$[D_0]$ of $\mathcal A$ with the $K^1({\mathcal A})$ class $[u]$.

Let $P$ denote the projection onto the nonnegative spectral subspace of $D_0$.
It is also well known 
 that spectral flow along $\{D_t^u\}$ is equal
to the Breuer-Fredholm index of the operator $PuP$
acting on $P\mathcal H$.
When ${\mathcal N}={\mathcal B}({\mathcal H})$
and the spectrum of $D_0$ is discrete \cite{CM}
show that $$ind(PuP)=\frac{1}{2}\tau_\omega(u^*[D_0,u]|D_0|^{-1}).$$
We aim to generalise this formula to 
the situation where $\mathcal N$ is a general semifinite
von Neumann algebra and link this formula with the expression for 
spectral flow.

\begin{lemma}
Let $D_0$ be an unbounded self-adjoint operator affiliated with $\mathcal N$
so that $(1+D_0^2)^{-1/2}$ is in
$\cl^{(1,\infty)}$. Let $A_t$ and $B$ be in $\mathcal N$ for $t\in[0,1]$
with $A_t$ self-adjoint and $t\mapsto A_t$ continuous. Let $D_t=D_0+A_t$ and
let $p$ be a real number with $1<p<4/3$.
Then, the quantity $$\tau\left (B(1+D_0^2)^{-p/2} - B(1+D_t^2)^{-p/2}\right )$$
is uniformly bounded independent of $t\in [0,1]$ and $p\in (1,4/3)$.
\end{lemma}

\begin{proof}
We estimate:
\begin{eqnarray*}
|\tau\left (B(1+D_0^2)^{-p/2} - B(1+D_t^2)^{-p/2}\right )| &\leq&
||B(1+D_0^2)^{-p/2} - B(1+D_t^2)^{-p/2}||_1 \\
&\leq& ||B||\cdot||(1+D_0^2)^{-p/2} - (1+D_t^2)^{-p/2}||_1 \\
&\leq& ||B||\cdot||\;|(1+D_0^2)^{-1} - (1+D_t^2)^{-1}|^{p/2}\;||_1 \\
&=& ||B||\cdot||(1+D_0^2)^{-1} - (1+D_t^2)^{-1}||^{p/2}_{p/2}.
\end{eqnarray*}
Where the last $inequality$ follows from the BKS inequality, see
\cite{BKS}, or the discussion and references in \cite{CPS}.

Now, by Lemma 2.9 of [CP1] we have
$$(1+D_0^2)^{-1} - (1+D_t^2)^{-1} = W_t + Z_t$$
where
$$W_t=D_0(1+D_0^2)^{-1}A_t(1+D_t^2)^{-1}$$
and
$$Z_t=(1+D_0^2)^{-1}A_tD_t(1+D_t^2)^{-1}.$$
Now, since $p/2$ is less than $1$, $||\cdot ||_{p/2}$ is not a norm: however,
by either
4.9 (iii) or 4.7 (i) of [FK] we have $$||W_t+Z_t||_{p/2}^{p/2} \leq
||W_t||_{p/2}^{p/2} + ||Z_t||_{p/2}^{p/2}.$$

Thus, it suffices to see that $||W_t||_{p/2}^{p/2}$ and $||Z_t||_{p/2}^{p/2}$
are bounded independent of $p$ and $t$.

Now, $(1+D_0^2)^{-1/2}$ and $(1+D_t^2)^{-1/2}$ are both in $\cl^{(1,\infty)}$
by Lemma 6 of [CP1] and therefore in $\cl^{q}$ for any $q>1$. In particular,
$(1+D_0^2)^{-1}$ and $(1+D_t^2)^{-1}$ are
both in $\cl^{2/3}$ and $\cl^{3/4}$.

Also, $p< 4/3$ implies $4-3p>0$ and since we also have $p>1$, we get
$r_p:=\frac{2p}{4-3p}>3/2$ and we easily calculate:
$$\frac{1}{2/3} + \frac{1}{r_p} = \frac{1}{p/2}.$$

So, by H\"older's inequality (Theorem 4.2 of [FK]), we get:
\begin{eqnarray*}
||W_t||_{p/2}^{p/2}&=&||D_0(1+D_0^2)^{-1}A_t(1+D_t^2)^{-1}||_{p/2}^{p/2}\\
&\leq& \left \{ ||D_0(1+D_0^2)^{-1}||_{r_p} ||A_t||\cdot
 ||(1+D_t^2)^{-1}||_{2/3}\right \}^{p/2}\\
&=& \left\{\left[\tau(|D_0(1+D_0^2)|^{-r_p})\right]^{1/r_p}||A_t||\left[\tau(
(1+D_t^2)^{-2/3}\right]^{3/2}\right\}^{p/2}\\
&\leq&\left\{\left[\tau((1+D_0^2)^{-r_p/2})\right]^{1/r_p}||A_t||f(||A_t||)
\left[\tau((1+D_0^2)^{-2/3}\right]^{3/2}\right\}^{p/2}
\end{eqnarray*}
where $f(t)=1+\frac{1}{2}(t^2+t\sqrt{4+t^2})$
by Lemma 6 of [CP1]. Since $r_p>3/2$ we have
$$\left[(1+D_0^2)^{-1/2}\right]^{3/2}
\geq\left[(1+D_0^2)^{-1/2}\right]^{r_p}.$$
Thus, we obtain our final inequality for $||W_t||_{p/2}^{p/2}$:
\begin{eqnarray*}
||W_t||_{p/2}^{p/2}&\leq&\left\{\left[\tau((1+D_0^2)^{-3/4})\right]^{1/r_p}||A_t||
f(||A_t||)\cdot||(1+D_0^2)^{-1}||_{2/3}\right\}^{p/2}.
\end{eqnarray*}
This last quantity is clearly a continuous function of $t$ and $p$ for
$t\in [0,1]$ and $p\in (1,4/3)$. As $p\to 1$ (and so $r_p\to 2$) we see that 
the estimate for $||W_t||_{p/2}^{p/2}$ converges to a continuous function
of $t\in [0,1]$ and so remains bounded at this end of $(1,4/3)$. On the other
hand, as $p\to 4/3$ (and so $r_p\to \infty$) we again see that the estimate 
for $||W_t||_{p/2}^{p/2}$ converges to a continuous function
of $t\in [0,1]$ and so remains bounded at the right hand end of $(1,4/3)$.
That is, the estimate for $||W_t||_{p/2}^{p/2}$ is bounded
independent of $t\in [0,1]$ and $p\in (1,4/3)$.

A slightly different calculation for $||Z_t||_{p/2}^{p/2}$, yields the
inequality:
$$||Z_t||_{p/2}^{p/2}\leq \left\{ ||(1+D_0^2)^{-1}||_{2/3}||A_t||\cdot
||f(A_t)||^{1/2}\left[||(1+D_0^2)^{-1}||_{3/4}^{3/4}\right]^{1/r_p}
\right\}^{p/2}.$$
Similar considerations to those above show that
$||Z_t||_{p/2}^{p/2}$ is also bounded independent of $t\in [0,1]$ and 
$p\in (1,4/3)$. This completes the proof.
\end{proof}

In (\cite{CP2} Corollary 9.4) we proved the following.
Let $\mathcal N$ be a factor and $({\mathcal N},D_0)$ be a
${\mathcal L}^{(1,\infty)}$-summable Breuer-Fredholm module for the unital
Banach $*$-algebra, $\mathcal A$, and let $u\in \mathcal A$
 be a unitary such that $[D_0,u]$
is bounded. Let $P$ be the projection on the non-negative spectral subspace of
$D_0$. Then for each $p>1$
$$
ind(PuP)=sf(D_0,uD_0u^*)= \frac{1}{\tilde C_{p/2}}\int_0^1\tau\left(u[D_0,u^*]
(1+(D_t^u)^2)^{-p/2}\right)dt$$
where
$$D_t^u = D_0 + tu[D_0,u^*] = D_0 + A_t\;\;\;\mbox{for}\;\;A_t=tu[D_0,u^*]\;\;
t\in[0,1]$$
and $\tilde C_{\frac{p}{2}}=\int_{-\infty}^{\infty} (1+x^2)^{-\frac{p}{2}}dx$.
(Note that a similar formula appears in Theorem 2.17 of \cite{CP1}
except that there the exponent $p>\frac{3}{2}$. The improvement
in the lower bound on the exponent uses the theory of theta summable
Fredholm modules in [CP2].)
The removal of the assumption that  $\mathcal N$ be a factor
is not hard (see for example the discussion in
the appendix to \cite{PR}).
The main point to note is that when  $\mathcal N$ is
a general semi-finite von Neumann algebra
then the map $u\to ind(PuP)$ is clearly
dependent on the choice of trace $\tau$, there being no
canonical choice. However this is not important for our
discussion in this paper.

\begin{thm}
Let $({\mathcal N},D_0)$ be a
${\mathcal L}^{(1,\infty)}$-summable Breuer-Fredholm module for the unital
Banach $*$-algebra, $\mathcal A$, and let $u\in \mathcal A$
 be a unitary such that $[D_0,u]$
is bounded. Let $P$ be the projection on the non-negative spectral subspace of
$D_0$. Then with $\omega$ chosen as in Theorem 1.5,
\begin{eqnarray*}
ind(PuP)=sf(D_0,uD_0u^*) &=& \lim_{p\to 1^+}\frac{1}{2}(p-1)\tau(u[D_0,u^*](1+D_0^2)^{-p/2})\\
& = & \frac{1}{2}\tau_{\omega}(u[D_0,u^*](1+D_0^2)^{-1/2})\\
& = & \frac{1}{2}\tau_{\omega}(u[D_0,u^*]|D_0|^{-1})
\end{eqnarray*}
where the last equality only holds if $D_0$ has a bounded inverse.
\end{thm}

\noindent{\bf Remarks}
(1) The equality
$$ind(PuP)=\frac{1}{2}\tau_{\omega}(u[D_0,u^*]|D_0|^{-1}) \eqno(6.1)$$
proved above should be compared with Theorem IV.2.8 of \cite{Co4}.
In the case where ${\mathcal N}={\mathcal B}({\mathcal H})$
the RHS of (6.1) is a Hochschild $1-$cocycle on $\mathcal A$
which is known to equal the Chern character of the
${\mathcal L}^{(1,\infty)}$-summable Fredholm module
$({\mathcal A},D_0,\mathcal H)$.

(2) Since any $1$-summable module is clearly a ${\mathcal L}^{(1,\infty)}$
-summable module, the theorem implies that any unbounded $1$-summable module
must have a trivial pairing with $K_1(\mathcal A)$ and is therefore 
uninteresting from the homological point of view.

\begin{proof}
By the extension of Corollary 9.4 of [CP2]
to the case where ${\mathcal N}$ is a general semifinite
von Neumann algebra, we have for each $p>1$, that
$$ind(PuP) = \frac{1}{\tilde C_{p/2}}\int_0^1\tau\left(u[D_0,u^*]
(1+(D_t^u)^2)^{-p/2}\right)dt$$
where the notation is described in the paragraph preceding the theorem.
Now, by Lemma 6.1, we have that
$$\left|\tau\left(u[D_0,u^*]
\left [(1+(D_t^u)^2)^{-p/2}-(1+D_0^2)^{-p/2}\right]\right)\right|$$
is uniformly bounded independent of $t$ and $p$ for $1<p<4/3$. Since,
$\tilde C_{p/2} \to\infty$ as $p\to 1^+$, we see that:
\begin{eqnarray*}
&&\left|ind(PuP) - \frac{1}{\tilde C_{p/2}}\tau\left(u[D_0,u^*](1+D_0^2)^{-p/2}
\right)\right|\\
&=&\left|\frac{1}{\tilde C_{p/2}}\int_0^1\tau\left(u[D_0,u^*]
(1+(D_t^u)^2)^{-p/2}\right)dt-
\frac{1}{\tilde C_{p/2}}\int_0^1\tau\left(u[D_0,u^*](1+D_0^2)^{-p/2}\right)
dt\right|\\
&\leq& \frac{1}{\tilde C_{p/2}}\int_0^1\left|\tau\left(u[D_0,u^*]
\left [(1+(D_t^u)^2)^{-p/2}-(1+D_0^2)^{-p/2}\right]\right)\right|dt\\
&\leq& \frac{Constant}{\tilde C_{p/2}} \to 0.
\end{eqnarray*}
Now, it is elementary that as $p\to 1^+$
$$\frac{2}{p-1}=\int_{|x|\geq 1}\left(\frac{1}{|x|}\right)^p dx \sim
\int_{-\infty}^{\infty}\left(\frac{1}{\sqrt{1+x^2}}\right)^p dx =
\tilde C_{p/2}.$$
This ends the proof of the first equality.

The second equality follows from Theorem 3.8(i).

The third equality follows from the fact that 
$\left(\sqrt{1+D_0^2}\right)^{-1} -|D_0|^{-1}$
is very trace-class:
$$\left(\sqrt{1+D_0^2}\right)^{-1} -|D_0|^{-1} =
\left(\sqrt{1+D_0^2}\right)^{-1} |D_0|^{-1}\left(\sqrt{1+D_0^2} +
|D_0|\right)^{-1}.$$
\end{proof}

\section{\bf Non-smooth foliations and pseudo-differential operators}

The main aim of Prinzis' thesis \cite{P} is to establish a
Wodzicki residue formula for the Dixmier trace
of certain pseudo-differential operators associated
to non-smooth actions of $\IR^n$ on a compact space $X$.
We will not reproduce the full details
of \cite{P}, indeed the subject deserves a far more
complete analysis than we have space for here.

The set-up is the group-measure space construction of Murray-von Neumann.
Thus $X$ is a compact space equipped with a probability measure $\nu$
and a continuous free minimal ergodic action $\alpha$ of $\IR^n$ on
$X$ leaving $\nu$ invariant. We write the action
as $x\to t.x$ for $x\in X$ and $t\in \IR^n$. Then the crossed product
$L^\infty(X,\nu)\times_\alpha\IR^n$ is a type $II$ factor
contained in the bounded operators on $L^2(\IR^n, L^2(X,\nu))$. We describe
the construction.
For a function $f\in L^1(\IR, L^\infty(X,\nu))
\subset L^\infty(X,\nu)\times_\alpha\IR^n$
the action of $f$ on a vector $\xi$
in $L^2(\IR^n, L^2(X,\nu))$ is defined
by twisted left convolution as follows:
$$(\tilde{\pi}(f)\xi)(s)
= \int_{\IR^n}\alpha_s^{-1}(f(t))\xi(s-t)dt.$$
Here $f(t)$ is a function on $X$ acting as a multiplication operator on
$ L^2(X,\nu)$.
The twisted convolution algebra
$$L^1(\IR^n, L^\infty(X,\nu)) \cap L^2(\IR^n, L^2(X,\nu))$$ is
a dense subspace of $L^2(\IR^n, L^2(X,\nu))$
and there is a canonical faithful, normal, semifinite trace, $Tr$,
on the von Neumann algebra that it generates. This von Neumann
algebra is $$\mathcal N=
(\tilde{\pi}(L^\infty(X,\nu)\times_{\alpha}\IR^n))^{\prime\prime}.$$
For functions $f,g: \IR^n\to L^\infty(X)$ which are in
$L^2(\IR^n, L^2(X,\nu))$ and whose twisted left convolutions $\tilde{\pi}(f),
\tilde{\pi}(g)$
define bounded operators on $L^2(\IR^n, L^2(X,\nu))$, this trace is given by:
$$Tr(\tilde{\pi}(f)^*\tilde{\pi}(g))=
\int_{\IR^n} \int_Xf(t,x)g(t,x)^*d\nu(x)dt$$
where we think of $f,g$ as functions on $\IR^n\times X$.

Identify
$L^2(\IR^n)$ with $L^2(\IR^n)\otimes 1\subset L^2(\IR, L^2(X,\nu))$
then any scalar-valued
function $f$ on $\IR^n$ which is the Fourier transform $f=\widehat{g}$ of a
bounded
$L^2$ function, $g$ will satisfy $f\in L^2(\IR, L^2(X,\nu))$
and $\tilde{\pi}(f)$ will be a bounded operator.

Pseudo-differential operators are defined in terms of their
symbols. A smooth symbol of order $m$ is a function
$a:X\times\IR^n\to\IC$ such that for each $x\in X$
$a_x$, defined by $a_x(t, \xi) = a(t.x, \xi)$,
satisfies

\noindent(1) $\sup\{|\partial_\xi^\beta\partial_t^\gamma a_x(t,\xi)(1+|\xi|)^{-m+|\beta|}
\ |\ (t,\xi)\in \IR^n\times\IR^n, \beta,\gamma\in \IN^n, |\beta|+|\gamma|
\leq M \}<\infty$
for all $M\in\IN$\\
(2)$\xi\to a_x(0,\xi)$ is a smooth function on $\IR^n$ into the space
${\mathcal C}^\infty(X)$, the set of continuous functions $f$ on $X$
such that $t \to(x\to f(t.x))$ is smooth on $\IR^n$.\\

Each symbol $a$ defines a
pseudo-differential operator $Op(a)$ on $C(X)\otimes C_c^\infty(\IR^n)$
by
$$Op(a)f(x,t)=\frac{1}{(2\pi)^n}\int_{\IR^n}e^{it\xi}
a(t.x,\xi)\hat f(x,\xi)d\xi,\ \ f\in  C(X)\otimes C_c^\infty(\IR^n).$$

The principal symbol of a pseudo-differential operator
$A$ on $X$ is the limit
$$\sigma_m(A)(x,\xi)=\lim_{\lambda\to\infty}\frac{a(x,\lambda\xi)}{\lambda^m}
(x,\xi)\in(X\times\IR^n\backslash\{0\})$$ if it exists.
We say $A$ is elliptic if its symbol $a$ is such that
$a_x$ is elliptic for all $x\in X$.
Prinzis studies invertible positive elliptic pseudo-differential
operators $A$ with a principal symbol. Henceforth we will
only consider such operators.
The zeta function of such an operator
is $\zeta(z)=\tau(A^z)$ and this exists because
$A^z$ is in the trace class in $\mathcal N$ \cite{P} for $\Re z<-n/m$.
Prinzis shows that
$$\lim_{x\to-\frac{n}{m}^-}(x+\frac{n}{m})\zeta(x)=
-\frac{1}{(2\pi)^nm}\int_{X\times S^{n-1}} \sigma_m(A)(x,\xi)^{-\frac{n}{m}}
d\nu(x)d\xi \eqno(7.1)$$
and that $A^{-\frac{n}{m}}\in {\mathcal L}^{(1,\infty)}$.

Our contribution to this situation is to note that (7.1)
combined with
Theorem 5.6 implies that we have the relation
$$\tau_\omega(A^{-\frac{n}{m}})=
\frac{1}{(2\pi)^nn}\int_{X\times S^{n-1}} \sigma_m(A)(x,\xi)^{-\frac{n}{m}}
d\nu(x)d\xi.$$
In other words we have a type $II$ Wodzicki residue
for evaluating the Dixmier trace of these
pseudo-differential operators.

\section{\bf Lesch's Index Theorem}
Here we consider a unital $C^*$-algebra $\mathcal A$
with a faithful finite trace,
$\tau$ satisfying $\tau(1)=1$ and a continuous action
 $\alpha$ of $\IR$ on $\mathcal A$
leaving $\tau$ invariant. In this section we deduce the index theorem of
M. Lesch as a corollary of our zeta function approach to the Dixmier Trace
formula for the index of generalised Toeplitz operators in this situation.
See \cite{L} and \cite{PR}.

We let $H_{\tau}$ denote the Hilbert space completion of $\mathcal A$
 in the inner
product $(a|b)=\tau(b^*a)$. Then
 $\mathcal A$ is a Hilbert Algebra and the left regular
representation of $\mathcal A$
 on itself extends by continuity to a representation,
$a \mapsto \pi_{\tau}(a)$ of $\mathcal A$
 on $H_{\tau}$ \cite{Dix}. In what follows, we will drop
the notation $\pi_{\tau}$ and just denote the action of $\mathcal A$
 on $H_{\tau}$
by juxtaposition.

We now look at the induced representation, $\tilde{\pi}$, of the crossed product
$C^*$-algebra
${\mathcal A}\times_{\alpha}\IR$ on $L^2(\IR,H_{\tau})$. That is, $\tilde{\pi}$ is the
representation $\pi\times\lambda$ obtained from the covariant pair,
$(\pi,\lambda)$ of
representations of the system $({\mathcal A},\IR,\alpha)$ defined
 for $a \in {\mathcal A}$,
$t,s \in \IR$ and $\xi \in L^2(\IR,H_{\tau})$ by:
$$(\pi(a)\xi)(s) = \alpha_s^{-1}(a)\xi(s)$$
and
$$\lambda_t(\xi)(s) = \xi(s-t).$$

Then, for a function $x\in L^1(\IR,{\mathcal A})
\subset {\mathcal A}\times_{\alpha}\IR$
 the action of $\tilde{\pi}(x)$ on a vector $\xi$
in $L^2(\IR,H_{\tau})$ is defined as follows:
$$(\tilde{\pi}(x)\xi)(s) = \int_{-\infty}^{\infty}\alpha_s^{-1}(x(t))\xi(s-t)dt.$$

Now the twisted convolution algebra
$L^1(\IR,{\mathcal A}) \cap L^2(\IR,H_{\tau})$ is
a dense subspace of $L^2(\IR,H_{\tau})$ and also a Hilbert Algebra in the
given inner product. As such, there is a canonical faithful, normal,
semifinite trace, $Tr$,
on the von Neumann algebra that it generates. Of course, this von Neumann
algebra is identical with
$$\mathcal N=
(\tilde{\pi}({\mathcal A}\times_{\alpha}\IR))^{\prime\prime}.$$
For functions $x,y: \IR\to {\mathcal A}\subset H_{\tau}$ which are in
$L^2(\IR,H_{\tau})$ and whose twisted left convolutions $\tilde{\pi}(x),
\tilde{\pi}(y)$
define bounded operators on $L^2(\IR,H_{\tau})$, this trace is given by:
$$Tr(\tilde{\pi}(y)^*\tilde{\pi}(x))= \langle x|y\rangle =
\int_{-\infty}^{\infty} \tau(x(t)y(t)^*)dt.$$

In particular, if we identify
$L^2(\IR)=L^2(\IR)\otimes 1_{\mathcal A}\subset L^2(\IR,H_{\tau})$
 then any scalar-valued
function $x$ on $\IR$ which is the Fourier transform $x=\widehat{f}$ of a
bounded
$L^2$ function, $f$ will have the properties that $x\in L^2(\IR,H_{\tau})$
and $\tilde{\pi}(x)$ is a bounded operator. For such scalar functions $x$, the
operator $\tilde{\pi}(x)$ is just the usual convolution by the function $x$ and is
usually denoted by $\lambda(x)$ since it is just the integrated form of
$\lambda$. The next Lemma follows easily from these considerations.

\begin{lemma}
With the hypotheses and notation discussed above\\

(i) if $h\in L^2(\IR)$ with $\lambda(h)$ bounded and $a\in \mathcal A$,
then defining $f:\IR\to H_{\tau}$ via $f(t)=ah(t)$ we see that
$f\in L^2(\IR,H_{\tau})$ and $\tilde{\pi}(f)=\pi(a)\lambda(h)$ is bounded,\\

(ii) if $g\in L^1(\IR)\cap L^{\infty}(\IR)$ and $a\in \mathcal A$ then
$\pi(a)\lambda(\hat{g})$ is trace-class in $\mathcal N$ and
$$Tr(\pi(a)\lambda(\hat{g}))=\tau(a)\int_{-\infty}^{\infty}g(t)dt.$$
\end{lemma}

\begin{proof}
To see part (i), let $\xi\in C_c(\IR,H_{\tau})\subseteq L^2(\IR,H_{\tau}).$
Then
\begin{eqnarray*}
(\tilde{\pi}(f)\xi)(s)&=&\int_{-\infty}^{\infty}\alpha_s^{-1}(f(t))\xi(s-t)dt\\
&=&\int_{-\infty}^{\infty}\alpha_s^{-1}(a)h(t)\xi(s-t)dt\\
&=&\alpha_s^{-1}(a)\int_{-\infty}^{\infty} h(t)\xi(s-t)dt\\
&=&\alpha_s^{-1}(a)(\lambda(h)\xi)(s)\\
&=&(\pi(a)\lambda(h)\xi)(s).
\end{eqnarray*}

To see part (ii) we can (and do) assume that $g$ is nonnegative and $a$ is
self-adjoint. Then let
$g=g^{1/2}g^{1/2}$ so that $g^{1/2}\in L^2\cap L^{\infty}$ and so
$\lambda(\widehat{g^{1/2}})$ is bounded. Now,
$$\pi(a)\lambda(\widehat{g}
)=\pi(a)\lambda(\widehat{g^{1/2}})\pi(1_{\mathcal A})\lambda
(\widehat{g^{1/2}}).$$
Then, $\pi(a)\lambda(\widehat{g^{1/2}})=\tilde{\pi}(x)$ where $x(t)=
a\widehat{g^{1/2}}(t)$
and $\pi(1_{\mathcal A})\lambda(\widehat{g^{1/2}})=\tilde{\pi}(y)$
 where $y(t)=
1_{\mathcal A}\widehat{g^{1/2}}(t).$
So, $\tilde{\pi}(x)$ and $\tilde{\pi}(y)$ are in ${\mathcal N}_{sa}$ and
$\pi(a)\lambda(\widehat{g})=\tilde{\pi}(x)\tilde{\pi}(y).$

Hence,
\begin{eqnarray*}
Tr(\pi(a)\lambda(\widehat{g}))&=&Tr(\tilde{\pi}(x)\tilde{\pi}(y))\\
&=&\int_{-\infty}^{\infty}\tau(x(t)y(t))dt\\
&=&\tau(a)\int_{-\infty}^{\infty}\left|\widehat{g^{1/2}}(t)\right|^2dt\\
&=&\tau(a)\int_{-\infty}^{\infty} g(s)ds.
\end{eqnarray*}

\end{proof}

Now, $\mathcal N$ is a semifinite von Neumann algebra with faithful, normal,
semifinite
trace, $Tr$, and a faithful representation $\pi:{\mathcal A}\to \mathcal N$
\cite{Dix}.
For each $t\in \IR$, $\lambda_t$ is a unitary in
$U(\mathcal N)$. In fact the one-parameter unitary group
$\{\lambda_t\; \vert\; t\in\IR\}$ can be written
$\lambda_t=e^{itD}$ where $D$  is
the unbounded
self-adjoint operator $$D=\frac{1}{2\pi i}\frac{d}{ds}$$ which is affiliated
with $\mathcal N$. In the Fourier Transform picture (i.e., the spectral picture for $D$)
of the previous proposition, $D$
becomes multiplication by the independent variable and so $f(D)$ becomes
pointwise multiplication by the function $f$. That is,
$$\tilde{\pi}(\hat{f}) = \lambda(\hat{f}) = f(D).$$
And, hence, if $f$ is a bounded $L^1$ function, then:
$$Tr(f(D)) = \int_{-\infty}^{\infty}f(t)dt.$$
By this discussion and the previous lemma, we have the following result

\begin{lemma}
If $f\in L^1(\IR)\cap L^{\infty}(\IR)$ and $a\in \mathcal A$
 then $\pi(a)f(D)$ is
trace-class in $\mathcal N$ and
$$Tr(\pi(a)f(D))=\tau(a)\int_{-\infty}^{\infty}f(t)dt.$$
\end{lemma}

 We let $\delta$ be the densely defined (unbounded) $*$-derivation
on $\mathcal A$ which is the infinitesimal generator of the representation
$\alpha:\IR\to Aut({\mathcal A})$
 and let $\hat{\delta}$ be the unbounded $*$-derivation
on $\mathcal N$ which is the infinitesimal generator of the
representation $Ad\circ\lambda:\IR\to Aut(\mathcal N)$
(here $Ad(\lambda_t)$ denotes conjugation by $\lambda_T$). Now if
$a\in dom(\delta)$ then clearly $\pi(a)\in dom(\hat{\delta})$ and
$\pi(\delta(a)) = \hat{\delta}(\pi(a))$. By \cite{BR} Proposition 3.2.55
(and its proof) we have that $\pi(\delta(a))$ leaves the domain of $D$
invariant and  $$\pi(\delta(a))=2\pi i[D,\pi(a)].$$
We are now in a position to state and prove Lesch's index theorem.

\begin{thm}
Let $\tau$ be a faithful finite trace on the unital $C^*$-algebra, $\mathcal A$,
which is invariant for an action $\alpha$ of $\IR$. Let $\mathcal N$ be the
semifinite von Neumann algebra
$(\tilde{\pi}({\mathcal A}\times_{\alpha}\IR))^{\prime\prime}$,
and let $D$ be the
infinitesimal generator of the canonical representation $\lambda$ of $\IR$
in $U(\mathcal N)$. Then, the representation
$\pi:{\mathcal A}\to \mathcal N$ defines a
${\mathcal L}^{(1,\infty)}$ summable Breuer-Fredholm module
$({\mathcal N},D)$ for
 $\mathcal A$. Moreover,
if $P$ is the nonnegative spectral projection for $D$ and
 $u\in U({\mathcal A})$ is also
in the domain of $\delta$, then $T_u:=P\pi(u)P$ is Breuer-Fredholm in
$P\mathcal NP$ and
$$ind (T_u) = \frac{1}{2\pi i}\tau(u\delta(u^*)).$$
\end{thm}

\begin{proof}
It is easy to see that $D$ satisfies
$(1+D^2)^{-1/2}\in {\mathcal L}^{(1,\infty)}$. By the previous discussion,
for any $a\in dom(\delta)$ we have $\pi(\delta(a))=2\pi i[D,\pi(a)].$ Since the
domain of $\delta$ is dense in $\mathcal A$ we see that $\pi$ defines a
${\mathcal L}^{(1,\infty)}$ summable Breuer-Fredholm module for $\mathcal A$.

Now, by Theorem 6.2 and Lemma 8.2
\begin{eqnarray*}
ind (T_u) &=& \lim_{p\to 1^+}\frac{1}{2}(p-1)
Tr(\pi(u)[D,\pi(u^*)](1+D^2)^{-p/2})\\
&=&\lim_{p\to 1^+}\frac{1}{2}(p-1)\frac{1}{2\pi i}
Tr(\pi(u\delta(u^*))(1+D^2)^{-p/2})\\
&=&\lim_{p\to 1^+}\frac{1}{2}(p-1)\frac{1}{2\pi i}
\tau(u\delta(u^*))\int_{-\infty}^{\infty}(1+t^2)^{-p/2}dt\\
&=&\lim_{p\to 1^+}\frac{1}{2\pi i}\tau(u\delta(u^*))\frac{1}{2}(p-1)
\tilde{C}_{p/2}\\
&=&\frac{1}{2\pi i}\tau(u\delta(u^*))
\end{eqnarray*}
\end{proof}


\begin{thebibliography}{APS3}

\bibitem [AS]{AS} M. Abramowitz and I.A. Stegun,
Handbook of mathematical functions with formulas, graphs, and
            mathematical tables, U.S. Govt. Print. Off, Washington, 1964.

\bibitem [BKS]{BKS} M.S. Birman, L.S. Koplienko and M.Z. Solomyak,
 \emph{Estimates
 for the spectrum of the difference between fractional powers of two
 selfadjoint operators,} Soviet Mathematics, {\bf 19}(3) (1975), 1-6.

\bibitem [BR]{BR} O. Bratteli and D. Robinson, Operator Algebras and Quantum
Statistical Mechanics I, Springer-Verlag, New York-Heidelberg-Berlin, 1979.

\bibitem[CM]{CM} A. Connes and H. Moscovici
\emph{The local index formula in non-commutative geometry},
Geometry and Functional Analysis {\bf 5} (1995) 174-243.

\bibitem[CP1]{CP1} A. L. Carey, J. Phillips, \emph{Unbounded Fredholm
Modules and Spectral Flow}, Canadian J. Math., {\bf 50}(4)(1998),
673--718.

\bibitem[CP2]{CP2} A. L. Carey, J. Phillips, \emph{Theta summable Fredholm
Modules, eta invariants and Spectral Flow}, Preprint.

\bibitem[CPS]{CPS} A. L. Carey, J. Phillips and F. A. Sukochev \emph{On
unbounded p-summable Fredholm modules}, Advances in Math., {\bf 151}
(2000) 140-163.

\bibitem[Co1]{Co1} A. Connes, Noncommutative Differential Geometry, Publ.
Math. Inst. Hautes Etudes Sci. (Paris), {\bfseries 62} (1985), 41--144.

\bibitem[Co2]{Co2} A. Connes, \emph{Cyclic Cohomology of Banach Algebras and
Characters of $\theta$-summable Fredholm Modules}, $K$-Theory {\bf1}(1988),
519--548.

\bibitem[Co3]{Co3} A. Connes, \emph{Compact Metric Spaces, Fredholm Modules
and Hyperfiniteness}, Ergodic Theory and Dynamical Systems {\bf9}(1989),
207--220.

\bibitem[Co4]{Co4} A. Connes, Noncommutative geometry, Academic Press,
San Diego, 1994.

\bibitem [CS]{CS} V.I.Chilin and F.A.Sukochev, {\it Weak convergence
in non-commutative symmetric spaces}, J. Operator Theory
{\bf 31} (1994), 35-65.

\bibitem[Dix]{Dix} J. Dixmier, Les alg\`ebres d'op\'erateurs dans l'espace
Hilbertien (Alg\`ebres de von Neumann), Gauthier-Villars, Paris, 1969.

\bibitem[Dix1]{Dix1} J. Dixmier,\emph{Existence de traces non-normales}
C.R. Acad. Sci. Paris {\bf A-B 262} (1966) A1107-A1108.

\bibitem[DPSS]{DPSS} P. G. Dodds, B. de Pagter, E. M. Semenov, F. A. Sukochev
\emph{Symmetric functionals and singular traces},
Positivity {\bf 2} (1998), 47-75.

\bibitem[E]{E} R.E. Edwards, Functional analysis: theory and applications,
Holt, Rinehart and Winston, New York, 1965.

\bibitem[F]{F} T. Fack, \emph{ Sur la notion de valeur caract\' eristique}, 
J. Operator Theory {\bf 7} (1982), 307--333.

\bibitem[FK]{FK} T. Fack and H. Kosaki, \emph{ Generalised $s$-numbers of
$\tau$-measurable operators}, Pacific J. Math. {\bf 123} (1986), 269--300.

\bibitem[G]{G} F. P. Greenleaf, Invariant Means on Topological Groups,
Van Nostrand Reinhold, New York, 1969.

\bibitem[GK]{GK} I. C. Gohberg, M. G. Krein, Introduction to the Theory
of Non-selfadjoint Operators, Translations of Mathematical Monographs,
{\bf vol. 18}, AMS, 1969.

\bibitem[Gr]{Gr} M. Gromov, \emph{K\"ahler-hyperbolicity and $L^2$
Hodge theory}, J. Diff. Geom.  {\bf  33}   (1991),   263-292.

\bibitem[GVF]{GVF} Gracia-Bondia, J.M., Varilly, J.C., Figueroa,
H. Elements of noncommutative geometry, Birkhauser, Boston, 2001.

\bibitem[H]{H} G.H. Hardy, Divergent Series, Clarendon, Oxford, 1949.

\bibitem[KR]{KR} R. Kadison and J. Ringrose,
  Fundamentals of the Theory of Operator Algebras,
Birkhauser, Boston, 1983.

\bibitem[L]{L} M. Lesch, \emph{On the Index of the Infinitesimal Generator of
a Flow}, J. Operator Theory {\bf 26} (1991), 73-92.

\bibitem[P]{P} R. Prinzis, \emph {Traces Residuelles
 et Asymptotique du Spectre d'Operateurs Pseudo-Differentiels}
Th\`ese, Universit\'e de Lyon, unpublished.

\bibitem[P1]{P1}J. Phillips,\emph {Self-adjoint Fredholm operators and spectral
flow}, Can. Math. Bull., {\bf 39}(1996), 460--467.

\bibitem[P2]{P2}J. Phillips,\emph {Spectral flow in type
I and II factors -- A new approach},
Fields Institute Communications, {\bf vol. 17}(1997), 137--153.

\bibitem[PR]{PR} J. Phillips and I. Raeburn, \emph{An Index Theorem for
Toeplitz Operators with Noncommutative Symbol Space}, J. Functional Analysis
{\bf 120} (1994), 239-263.

\end{thebibliography}
\end{document}